\documentclass[reqno]{amsart}

\usepackage{amsmath}
\usepackage{amssymb}
\usepackage{amsfonts}
\usepackage{amsthm}
\usepackage[foot]{amsaddr}
\usepackage{mathtools}
\mathtoolsset{%
above-intertext-sep = -1ex
}

\usepackage[utf8]{inputenc}
\usepackage[T1]{fontenc}

\usepackage[sf,mono=false]{libertine}
\usepackage[sans]{dsfont}
\usepackage[%
cal=cm,
scr=euler,
]
{mathalfa}

\usepackage[dvipsnames,svgnames]{xcolor}
\colorlet{MyBlue}{DodgerBlue!75!Black}
\colorlet{MyGreen}{DarkGreen!95!Black}

\usepackage[font=small,labelfont=bf]{caption}
\usepackage{subfigure}

\usepackage{acronym}
\usepackage{latexsym}
\usepackage{paralist}
\usepackage{wasysym}
\usepackage{xspace}

\usepackage[numbers,sort&compress]{natbib}

\usepackage{hyperref}
\hypersetup{
colorlinks=true,
linktocpage=true,
pdfstartview=FitH,
breaklinks=true,
pdfpagemode=UseNone,
pageanchor=true,
pdfpagemode=UseOutlines,
plainpages=false,
bookmarksnumbered,
bookmarksopen=false,
bookmarksopenlevel=1,
hypertexnames=true,
pdfhighlight=/O,
urlcolor=MyBlue!60!black,linkcolor=MyBlue!70!black,citecolor=DarkGreen!70!black, % <--- for screen
pdftitle={},
pdfauthor={},
pdfsubject={},
pdfkeywords={},
pdfcreator={pdfLaTeX},
pdfproducer={LaTeX with hyperref}
}
\newcommand{\EMAIL}[1]{\email{\href{mailto:#1}{#1}}}

\numberwithin{equation}{section}  %numberwithin goes before cleverefs when using hyperref
\usepackage[sort&compress,capitalize,nameinlink]{cleveref}
\crefrangeformat{equation}{\upshape(#3#1#4)\textendash(#5#2#6)}

\newcommand{\dd}{\:d}

\newcommand{\eps}{\varepsilon}
\newcommand{\from}{\colon}

\newcommand{\dif}{\dd}

\newcommand{\R}{\mathbb{R}}

\newcommand{\N}{\mathbb{N}}

\DeclareMathOperator*{\argmin}{argmin}
\DeclareMathOperator*{\argmax}{argmax}
 
\DeclareMathOperator{\bd}{bd}
\DeclareMathOperator{\cl}{cl}

\DeclareMathOperator{\diam}{diam}

\DeclareMathOperator{\dom}{dom}

\DeclareMathOperator{\Int}{ri}

\DeclareMathOperator{\image}{im}

\DeclareMathOperator{\tr}{tr}

\DeclareMathOperator{\Id}{Id}

\newcommand{\Leb}{{\mathsf{Leb}}}

\newcommand{\bT}{\mathbf{T}}

\renewcommand{\emptyset}{\varnothing}

\newcommand{\Rn}{\R^n}
\newcommand{\Rnp}{\R_{+}^{n}}
\newcommand{\Rp}{\R_+}

\renewcommand{\Pr}{\mathbb{P}}

\newcommand{\scrD}{\mathcal{D}}

\newcommand{\scrF}{\mathcal{F}}

\newcommand{\scrL}{\mathcal{L}}

\newcommand{\scrO}{\mathcal{O}}

\newcommand{\scrQ}{\mathcal{Q}}

\newcommand{\scrV}{\mathcal{V}}

\newcommand{\scrX}{\mathcal{X}}
\newcommand{\scrY}{\mathcal{Y}}

\theoremstyle{plain}
\newtheorem{theorem}{Theorem}
\newtheorem{corollary}[theorem]{Corollary}
\newtheorem*{corollary*}{Corollary}
\newtheorem{lemma}[theorem]{Lemma}
\newtheorem{proposition}[theorem]{Proposition}

\theoremstyle{definition}
\newtheorem{definition}[theorem]{Definition}
\newtheorem*{definition*}{Definition}

\renewcommand\qed{\hfill\small$\blacksquare$}
\theoremstyle{remark}
\newtheorem{remark}{Remark}
\newtheorem*{remark*}{Remark}
\newtheorem*{notation*}{Notational remark}
\newtheorem{example}{Example}

\numberwithin{theorem}{section}
\numberwithin{remark}{section}
\numberwithin{example}{section}

\newcommand{\as}{\textup(a.s.\textup)\xspace}
\newcommand{\textpar}[1]{\textup(#1\textup)}

\DeclareMathOperator{\ex}{\mathbb{E}}
\DeclareMathOperator{\prob}{\mathbb{P}}
\DeclareMathOperator{\VI}{VI}
\DeclareMathOperator{\val}{val}

\DeclarePairedDelimiter{\braces}{\{}{\}}
\DeclarePairedDelimiter{\bracks}{[}{]}
\DeclarePairedDelimiter{\parens}{(}{)}

\DeclarePairedDelimiter{\abs}{\lvert}{\rvert}

\DeclarePairedDelimiter{\inner}{\langle}{\rangle}
\DeclarePairedDelimiter{\norm}{\lVert}{\rVert}
\DeclarePairedDelimiter{\dnorm}{\lVert}{\rVert_{\ast}}

\DeclarePairedDelimiterX{\braket}[2]{\langle}{\rangle}{#1,#2}

\DeclarePairedDelimiterX{\setdef}[2]{\{}{\}}{#1:#2}

\DeclarePairedDelimiterXPP{\probof}[1]{\prob}{(}{)}{}{%

#1}

\DeclarePairedDelimiterXPP{\exof}[1]{\ex}{[}{]}{}{%

#1}

\newcommand{\cost}{c}
\newcommand{\sol}{x_{\ast}}
\newcommand{\olim}{\hat x}
\newcommand{\bigoh}{\scrO}

\newcommand{\noisevol}{\sigma}
\newcommand{\noisedev}{\noisevol_{\ast}}
\newcommand{\noisevar}{\noisedev^{2}}

\newcommand{\nhd}{\scrO}

\begin{document}

%*************************************************************
%*****    FRONT MATTER AND METADATA
%*************************************************************

%----------------------------------------------------------------------
%%% TITLE & AUTHORS
%----------------------------------------------------------------------
\title[Stochastic mirror descent in variational inequalities]{Stochastic mirror descent dynamics and their convergence in monotone variational inequalities}

\author{Panayotis Mertikopoulos$^{\ast}$}
\address
{$^{\ast}$\,%
Univ. Grenoble Alpes, CNRS, Inria, LIG, F-38000, Grenoble, France.}
\EMAIL{panayotis.mertikopoulos@imag.fr}
%\URLADDR{http://mescal.imag.fr/membres/panayotis.mertikopoulos}

\author{Mathias Staudigl$^{\sharp}$}
\address{$^{\sharp}$\,%
Maastricht University, Department of Quantitative Economics, P.O. Box 616, NL\textendash 6200 MD Maastricht, The Netherlands.}
\EMAIL{m.staudigl@maastrichtuniversity.nl}

%----------------------------------------------------------------------
%%% KEYWORDS
%----------------------------------------------------------------------

%----------------------------------------------------------------------
%%% ACRONYMS
%----------------------------------------------------------------------
\newcommand{\acli}[1]{\textit{\acl{#1}}}
\newcommand{\acdef}[1]{\textit{\acl{#1}} \textup{(\acs{#1})}\acused{#1}}
\newcommand{\acdefp}[1]{\emph{\aclp{#1}} \textup(\acsp{#1}\textup)\acused{#1}}

\newacro{APT}{asymptotic pseudotrajectory}
\newacroplural{APT}[APTs]{asymptotic pseudotrajectories}
\newacro{NE}{Nash equilibrium}
\newacroplural{NE}[NE]{Nash equilibria}
\newacro{MD}{mirror descent}
\newacro{SMD}{stochastic mirror descent}
\newacro{SDE}{stochastic differential equation}
\newacro{VI}{variational inequality}
\newacroplural{VI}[VIs]{variational inequalities}
\newacro{iid}[i.i.d.]{independent and identically distributed}

\begin{abstract}
We examine a class of stochastic mirror descent dynamics in the context of monotone variational inequalities (including Nash equilibrium and saddle-point problems).
The dynamics under study are formulated as a stochastic differential equation, driven by a (single-valued) monotone operator and perturbed by a Brownian motion. The system's controllable parameters are two variable weight sequences, that respectively pre- and post-multiply the driver of the process. By carefully tuning these parameters, we obtain global convergence in the ergodic sense, and we estimate the average rate of convergence of the process. We also establish a large deviations principle, showing that individual trajectories exhibit exponential concentration around this average.
\end{abstract}
\maketitle
\section{Introduction}
\label{sec:introduction}
Dynamical systems governed by monotone operators play an important role in the fields of optimization (convex programming), game theory (Nash equilibrium and generalized Nash equilibrium problems), fixed-point theory, partial differential equations, and many other fields of applied mathematics. In particular, the study of the relationship between continuous- and discrete-time models has given rise to a vigorous literature at the interface of these fields -  see, e.g., \cite{PeySor10} for a recent overview and \cite{WibWilJor16} for connections to accelerated methods.\\
The starting point of much of this literature is that an iterative algorithm can be seen as a discretization of a continuous dynamical system. Doing so sheds new light on the properties of the algorithm, it provides Lyapunov functions, which are useful for its asymptotic analysis, and often leads to new classes of algorithms altogether. A classical example of this arises in the study of (projected) gradient descent dynamics and its connection with Cauchy's steepest descent algorithm \textendash\ or, more generally, in the relation between the mirror descent (MD) class of algorithms \cite{NY83} and dynamical systems derived from Bregman projections and Hessian Riemannian metrics \cite{BT03,ABB04,ABRT04}.

\section{Problem Formulation and Related Literature}

Throughout this paper, $\scrX$ denotes a compact convex subset of an $n$-dimensional real space $\scrV\cong\Rn$ with norm $\norm{\cdot}$. We will also write $\scrY\equiv\scrV^{\ast}$ for the dual of $\scrV$,
$\braket{y}{x}$ for the canonical pairing between $y\in\scrV^{\ast}$ and $x\in\scrV$,
and $\dnorm{y}:=\sup\setdef{\braket{y}{x}}{\norm{x} \leq 1}$ for the dual norm of $y$ in $\scrV^{\ast}$.
We denote the relative interior of $\scrX$ by $\Int(\scrX)$, and its boundary by $\bd(\scrX)$.\\
In this paper we are interested in deriving dynamical system approaches to solve monotone variational inequality (VI) problems. To define them, let $v:\scrX\to\scrY$ be a Lipschitz continuous monotone map, i.e.
 \begin{itemize}
 \label{eq:Lip}
 \item[(H1)] $ \dnorm{v(x)-v(x')}
	\leq L \norm{x-x'},\text{ and }\inner{v(x)-v(x'),x-x'}\geq 0,$ for some $L>0$ and for all $x,x'\in\scrX$.
	\end{itemize}
Throughout this paper, we will be interested in solving the Minty VI problem:
\begin{equation}
\label{eq:VI}
\tag{MVI}
\text{Find $\sol\in\scrX$ such that $\inner{v(x),x-\sol} \geq 0$ for all $x\in\scrX$}.
\end{equation}
Since $v$ is assumed continuous and monotone, this VI problem is equivalent to the \emph{Stampacchia} VI problem:\footnote{Equivalence is understood in the sense that $x_{\ast}\in\scrX$ solves \eqref{eq:VI} if and only if it solves \eqref{eq:VI-strong}. See Proposition 3.1 in \cite{HarPan90}, and Proposition 3 in \cite{Gia1998}, which treats also extensions to vector Variational Inequality problems.}
\begin{equation}
\label{eq:VI-strong}
\tag{SVI}
\text{Find $\sol\in\scrX$ such that $\inner{v(\sol),x-\sol} \geq 0$ for all $x\in\scrX$}.
%\text{Find $\sol\in\scrX$ s.t. } \inner{v(\sol),x-\sol}\geq 0\quad\forall x\in\scrX.
\end{equation}
When we need to keep track of $\scrX$ and $v$ explicitly, we will refer to \eqref{eq:VI} and/or \eqref{eq:VI-strong} as $\VI(\scrX,v)$.
The solution set of $\VI(\scrX,v)$ will be denoted as $\scrX_{\ast}$;
by standard results, $\scrX_{\ast}$ is convex, compact and nonempty \cite{FacPan03}. Below, we present a selected sample of examples and applications of VI problems;
for a more extensive discussion, see \cite{FerPan97,FacPan03,ScuPalFacPan10}.

\begin{example}[Convex optimization]
\label{ex:opt}
Consider the problem
\begin{equation}
\label{eq:opt}
\tag{Opt}
\begin{aligned}
\textrm{minimize}
	\quad
	&f(x),
	\\
\textrm{subject to}
	\quad
	&x\in\scrX,
\end{aligned}
\end{equation}
where $f\from\scrX\to\R$ is convex and continuously differentiable on $\scrX$.
If $\sol$ is a solution of \eqref{eq:opt}, first-order optimality gives
\begin{equation}
\inner{\nabla f(\sol),x-\sol}
	\geq 0
	\quad
	\text{for all $x\in\scrX$}.
\end{equation}
Since $f$ is convex, $v= \nabla f$ is monotone, so \eqref{eq:opt} is equivalent to $\VI(\scrX,\nabla f)$ \cite{RW98}. 
\end{example}
\medskip

\begin{example}[Saddle-point problems]
\label{ex:saddle}
Let $\scrX^{1}\subseteq\R^{n_{1}}$ and $\scrX^{2}\subseteq\R^{n_{2}}$ be compact and convex, and let $U\from\scrX^{1}\times\scrX^{2}\to\R$ be a smooth convex-concave function (i.e. $U(x^{1},x^{2})$ is convex in $x^{1}$ and concave in $x^{2}$).
Then, the associated \emph{saddle-point} (or \emph{min-max}) problem is to determine the \emph{value} of $U$, defined here as
%the unique number $\val$ satisfying 
\begin{equation}
\label{eq:value}
\tag{Val}
\val
	= \min_{x^{1}\in\scrX^{1}} \max_{x^{2}\in\scrX^{2}} U(x^{1},x^{2})
	= \max_{x^{2}\in\scrX^{2}} \min_{x^{1}\in\scrX^{1}} U(x^{1},x^{2}).
\end{equation}
Existence of $\val$ follows directly from von Neumann's minimax theorem. Moreover, letting
\begin{equation}
\label{eq:v-saddle}
v(x^{1},x^{2}):= \big(\nabla_{x^{1}} U(x^{1},x^{2}), - \nabla_{x^{2}} U(x^{1},x^{2}) \big),
\end{equation}
it is easy to check that $v$ is monotone as a map from $\scrX:=\scrX^{1}\times\scrX^{2}$ to $\R^{n_{1}+n_{2}}$ (because $U$ is convex in its first argument and concave in the second). Then, as in the case of \eqref{eq:opt}, first-order optimality implies that the saddle-points of \eqref{eq:value} are precisely the solutions of $\VI(\scrX,v)$ \cite{Nes09}.
\end{example}
\medskip

\begin{example}[Convex games]
\label{ex:games}
One of the main motivations for this paper comes from determining Nash equilibria of games with convex cost functions.
To state the problem, let $\mathcal{N} := \{1,\dotsc,N\}$ be a finite set of \emph{players} and, for each $i\in\mathcal{N}$, let $\scrX^{i}\subseteq\R^{n_{i}}$ be a compact convex set of \emph{actions} that can be taken by player $i$.
Given an action profile $x = (x^{1},\dotsc,x^{N}) \in\scrX:=\prod_{i} \scrX^{i}$, the cost for each player is determined by an associated \emph{cost function} $\cost^{i}\from\scrX\to\R$.
The unilateral minimization of this cost leads to the notion of Nash equilibrium, defined here as an action profile $\sol = (\sol^{i})_{i\in\mathcal{N}}$ such that
\begin{equation}
\label{eq:Nash}
\tag{NE}
\cost^{i}(\sol)
	\leq \cost^{i}(x^{i};\sol^{-i})
	\quad
	\text{for all $x^{i}\in\scrX^{i}$, $i\in\mathcal{N}$}.
\end{equation}
Of particular interest to us is the case where each $\cost^{i}$ is smooth and individually convex in $x^{i}$.
In this case, the profile $v(x) = (v^{i}(x))_{i\in\mathcal{N}}$ of individual gradients $v^{i}(x):= \nabla_{x^{i}} \cost^{i}(x)$ forms a monotone map and, by first-order optimality, the Nash equilibrium problem \eqref{eq:Nash} boils down to solving $\VI(\scrX,v)$ \cite{FacPan03,Mer16}.
\end{example}
\medskip

In the rest of this paper, we will consider two important special cases of operators $v:\scrX\to\scrV^{\ast}$, namely:
\begin{enumerate}
\item
\emph{Strictly monotone} problems, i.e. when
\begin{alignat}{3}
\label{eq:strict}
\inner{v(x') - v(x),x'-x}
	&\geq 0
	&\quad
	&\text{with equality if and only if $x=x'$}.
\intertext{%
\item
\emph{Strongly monotone} problems, i.e. when}
\label{eq:strong}
\inner{v(x')-v(x),x'-x}
	&\geq \gamma\norm{x'-x}^{2}
	&\quad
	&\text{for some $\gamma>0$}. 
\end{alignat}
\end{enumerate}
Clearly, strong monotonicity implies strict monotonicity (which in turns implies ordinary monotonicity).
In the case of convex optimization problems, strict (respectively strong) monotonicity corresponds to strict (respectively strong) convexity of the problem's objective function. Under either refinement, \eqref{eq:VI} admits a unique solution, which will be referred to as ``the'' solution of \eqref{eq:VI}.

\subsection{Contributions}
Building on the above, this paper is concerned with a \emph{stochastic} dynamical system resulting from Nesterov's well-known ``dual-averaging'' mirror descent algorithm \cite{Nes09} perturbed by noise and/or random disturbances. Heuristically, this algorithm aggregates descent steps in the problem's (unconstrained) dual space, and then ``mirrors'' the result back to the problem's feasible region to obtain a candidate solution at each iteration.
This ``mirror step'' is performed as in the classical setting of \cite{NY83,BT03}, but the dual aggregate is further post-multiplied by a variable parameter (thus turning ``dual aggregates'' into ``dual averages'').  Thanks to this averaging, the resulting algorithm is particularly suited for problems where only noisy information is available to the optimizer, rendering it particularly useful for machine learning and engineering applications \cite{Xia10}, even when the stochastic environment is not stationary \cite{DucAgaWai12}.\footnote{In the deterministic case, similar results are to be contrasted with the use of Gauss-Seidel methods for solving Variational and (generalized) Nash Equilibrium problems; for a survey, see \cite{FK07,SFPP10}.}\\
In more detail, the dynamics under study are formulated as a stochastic differential equation (SDE) driven by a (single-valued) monotone operator and perturbed by an It\^{o} martingale noise process.
As in Nesterov's original method \cite{Nes09}, the dynamics' controllable parameters are two variable weight sequences that respectively pre- and post-multiply the drift of the process:
the first acts as a ``step-size'' of sorts, whereas the second can be seen as an ``inverse temperature'' parameter (as in simulated annealing).
By carefully tuning these parameters, we are then able to establish the following results:
First, if the intensity of the noise process decays with time, the dynamics converge to the (deterministic) solution of the underlying VI problem (cf.~\cref{sec:smallnoise}).
Second, in the spirit of the ergodic convergence analysis of \cite{Nes09}, we establish that this convergence can be achieved at an $\bigoh(1/\sqrt{t})$ rate on average (\cref{sec:ergodic}).%
\footnote{Interestingly, the corresponding rate for the deterministic (noise-free) dynamics is $\bigoh(1/t)$, indicating a substantial drop from the deterministic to the stochastic regime.
This drop is consistent with the black-box convergence rate of Mirror Descent in (stochastic) VI problems \cite{Nes09} and is due to the second-order It\^{o} correction that appears in the stochastic case.}
Finally, in \cref{sec:LDP}, we establish a large deviations principle showing that, as far as ergodic convergence is concerned, the above convergence rate holds with (exponentially) high probability, not only in the mean.\\
Conceptually, our work here has close ties to the literature on dynamical systems that arise in the solution of VI problems, see e.g. \cite{BT03,CEG09,GP14,AA15,KBB15,WibWilJor16}, and references therein.
More specifically, a preliminary version of the dynamics considered in this paper was recently studied in the context of convex programming and gradient-like flows in \cite{RB13,MerStaConvex}.
The ergodic part of our analysis here extends the results of \cite{MerStaConvex} to saddle-point problems and monotone variational inequalities, while the use of two variable weight sequences allows us to obtain almost sure convergence results without needing to rely on a parallel-sampling mechanism for variance reduction as in \cite{RB13}.

\subsection{Stochastic Mirror Descent Dynamics}
\label{sec:SMD}
Mirror descent is an iterative optimization algorithm combining first-order oracle steps with a ``mirror step'' generated by a projection-type mapping. For the origins of the method, see \cite{NY83}. The key ingredient defining this mirror step is a generalization of the Euclidean distance known as a ``distance-generating'' function:
\begin{definition}
\label{def:dgf}
We say that $h\from\scrX\to\R$ is a \emph{distance-generating function} on $\scrX$, if
\begin{enumerate}
[\itshape a\upshape)]
\item
$h$ is \emph{continuous}.
\item
$h$ is \emph{strongly convex}, i.e. there exists some $\alpha>0$ such that
\begin{equation}
\label{eq:StrongConvex}
h(\lambda x+(1-\lambda)x')
	\leq \lambda h(x)+(1-\lambda)h(x')-\frac{\alpha}{2}\lambda(1-\lambda)\norm{x-x'}^{2},
\end{equation}
for all $x,x'\in\scrX$ and all $\lambda\in[0,1]$.
\end{enumerate}
\end{definition}

Given a distance-generating function on $\scrX$, its convex conjugate is given by
\begin{flalign}
h^{\ast}(y):=
	\max_{x\in\scrX} \{\inner{y,x} - h(x) \},
	\quad
	y\in\scrY,
\end{flalign}
and the induced \emph{mirror map} is defined as
\begin{equation}
\label{eq:Q}
Q(y):= \argmax_{x\in\scrX} \{ \inner{y,x} - h(x) \}.
\end{equation}
Thanks to the strong convexity of $h$, $Q(y)$ is well-defined and single-valued for all $y\in\scrY$.
In particular, as illustrated in the examples below, it plays a role similar to that of a projection mapping:
\medskip

\begin{example}[Euclidean distance]
If $h(x) = \frac{1}{2} \norm{x}_{2}^{2}$, the induced mirror map is the standard Euclidean projector
\begin{flalign}
%\label{eq:Eucl}
Q(y)
	= \argmax_{x\in\scrX} \braces[\bigg]{\sum\nolimits_{j=1}^{n} y_{j} x_{j} - \frac{1}{2} \sum\nolimits_{j=1}^{n} x_{j}^{2} }
	= \argmin_{x\in\scrX} \norm{x-y}_{2}^{2}.
\end{flalign}
\end{example}
\medskip

\begin{example}[Gibbs\textendash Shannon entropy]
If $\scrX = \setdef{x\in\Rnp}{\sum_{j=1}^{n} x_{j}=1}$ is the unit simplex in $\Rn$, the (negative) Gibbs-Shannon entropy $h(x) = \sum_{j=1}^{n} x_{j} \log x_{j}$ gives rise to the so-called \emph{logit choice} map
\begin{flalign}
%\label{eq:logit}
Q(y)
	= \frac{(\exp(y_{j}))_{j=1}^{n}}{\sum_{k=1}^{n} \exp(y_{k})}.
\end{flalign}
\end{example}
\medskip

\begin{example}[Fermi\textendash Dirac entropy]
If $\scrX = [0,1]^{n}$ is the unit cube in $\Rn$, the (negative) Fermi-Dirac entropy $h(x) = \sum_{j=1}^{n} [x_{j} \log(x_{j}) + (1-x_{j})\log(1-x_{j})]$ induces the so-called \emph{logistic map}
\begin{flalign}
%\label{eq:FermiDirac}
Q(y)
	= \parens*{\frac{\exp(y_{j})}{1 + \exp(y_{j})}}_{j=1}^{n}
\end{flalign}
\end{example}
\medskip

For future reference, some basic properties of mirror maps are collected below:
\begin{proposition}
\label{prop:Q}
Let $h$ be a distance-generating function on $\scrX$.
Then, the induced mirror map $Q\from\scrY\to\scrX$ satisfies the following properties:
\begin{enumerate}
[\indent\itshape a\upshape)]
\item
$x=Q(y)$ if and only if $y\in\partial h(x)$, where 
\begin{flalign}
\partial h(x):=\{p\in\scrV^{\ast} : h(y)\geq h(x)+\inner{p,y-x}\quad\forall y\in\scrX\},
\end{flalign}
is the subgradient of $h$ at $x$. In particular, $\image Q=\dom\partial h=\{x\in\scrX : \partial h(x)\neq\emptyset\}$.
\item
$h^{\ast}$ is continuously differentiable on $\scrY$ and $\nabla h^{\ast}(y) = Q(y)$ for all $y\in\scrY$.
\item $Q(\cdot)$ is $(1/\alpha)$-Lipschitz continuous.
\end{enumerate}
\end{proposition}

The properties reported above are fairly standard in convex analysis;
for a proof, see e.g. \cite[p.~217]{Zal02}, \cite[Theorem~23.5]{Roc70} and \cite[Theorem~12.60(b)]{RW98}.
Of particular importance is the minimizing argument identity $\nabla h^{\ast} = Q$ which provides a quick way of calculating $Q$ in ``prox-friendly'' geometries (such as the examples discussed above).\\
Now, as mentioned above, mirror descent exploits the flexibility provided by a (not necessarily Euclidean) mirror map by using it to generate first-order steps along $v$.
For concreteness, we will focus on the so-called ``dual averaging'' variant of mirror descent \cite{Nes09}, defined here via the recursion
\begin{equation}
\label{eq:MD-discrete}
y_{t+1}
	= y_{t} - \lambda_{t} v(x_{t})
	\qquad 
x_{t+1}
	= Q(\eta_{t+1}y_{t+1}),
\end{equation}
where:
\begin{enumerate}
[\indent 1)]
\item
$t=0,1,\dotsc$ denotes the stage of the process.
\item
$y_{t}$ is an auxiliary dual variable that aggregates first-order steps along $v$.%
\footnote{The usual initialization is $y_{0}=0$, $x_{0} = Q(0) = \argmin_{x\in\scrX} h(x)$, but other initializations are possible.}
\item
$\lambda_{t}$ is a variable step-size parameter that \emph{pre-multiplies} the input at each stage.
\item
$\eta_{t}$ is a variable weight parameter that \emph{post-multiplies} the dual aggregate $y_{t}$.%
\footnote{The name ``dual averaging'' alludes to the choice $\lambda_{t} = 1$, $\eta_{t} = 1/t$:
under this choice of parameters, $x_{t}$ is a mirror projection of the ``dual average'' $y_{t} = t^{-1} \sum_{s=0}^{t-1} v(x_{s})$.}
\end{enumerate}

Passing to continuous time, we obtain the \emph{mirror descent dynamics}
\begin{equation}
\label{eq:MD}
\tag{MD}
dy(t)
	= -\lambda(t)\, v(x(t)) \dif t
	\qquad
x(t)
	= Q(\eta(t) y(t)),
\end{equation}
with $\eta(t)$ and $\lambda(t)$ serving the same role as before (but now defined over all $t\geq0$).
In particular, our standing assumption for the parameters $\lambda$ and $\eta$ of \eqref{eq:MD} will be that
\begin{itemize}
\label{eq:params}
\item[(H2)] $\eta(t)$ and $\lambda(t)$ are positive, $C^{1}$-smooth and nonincreasing.
\end{itemize}
Heuristically, the assumptions above guarantee that the dual process $y(t)$ does not grow too large too fast, so blow-ups in finite time are not possible.
%\footnote{Note that we are not assuming that $\lambda(t)$ is nonincreasing in $t$ (as is typically the case for step-size parameters).
%The reason for this is that $\eta(t)$ provides an effective counter-balance to $\lambda(t)$ which can be exploited to yield accelerated rates of convergence;
%for a detailed discussion, see \cref{sec:stochastic} below.}
Together with the basic convergence properties of the dynamics \eqref{eq:MD}, this is discussed in more detail in \cref{sec:deterministic} below.
%The analysis and basic convergence properties of the dynamics \eqref{eq:MD} are presented in \cref{sec:deterministic} below.

The primary case of interest in our paper is when the oracle information for $v(x)$ in \eqref{eq:MD} is subject to noise, measurement errors and/or other stochastic disturbances.
To account for such perturbations, we will instead focus on the \emph{stochastic mirror descent dynamics}
\begin{equation}
\label{eq:SMD}
\tag{SMD}
dY(t)
	= -\lambda(t)\, \bracks{v(X(t)) \dif t + \dif M(t)}
	\qquad
X(t)
	= Q(\eta(t)Y(t))
\end{equation}
where $M(t)$ is a continuous martingale with respect to some underlying stochastic basis $(\Omega,\scrF,(\scrF_{t})_{t\geq 0},\Pr)$.%
\footnote{We tacitly assume here that the filtration $(\scrF_{t})_{t\geq 0}$ satisfies the usual conditions of right continuity and completeness, and carries a standard $d$-dimensional Wiener process $(W(t))_{t\geq 0}$.}
In more detail, we assume for concreteness that the stochastic disturbance term $M(t)$ is an It\^{o} process of the form
\begin{equation}
\label{eq:noise}
dM(t)
	= \noisevol(X(t),t) \cdot \dif W(t),
\end{equation}
where $W(t)$ is a $d$-dimensional Wiener process adapted to $\scrF_{t}$, and $\noisevol(x,t)$ is an $n\times d$ matrix capturing the \emph{volatility} of the noise process. Heuristically, the volatility matrix of $M(t)$ captures the intensity of the noise process and the possible correlations between its components.
%For instance, when $d=n$ and $\noisevol$ is the identity matrix, $M(t)$ is just a standard Wiener process:
%in this case, the increments of the noise are \ac{iid} and they are not correlated across different components.
%Otherwise, if $\noisevol$ is not diagonal, $M(t)$ could exhibit nontrivial correlations and/or other dependencies across components.

In terms of regularity, we will be assuming throughout that $\noisevol(x,t)$ is measurable in $t$, as well as bounded, and Lipschitz continuous in $x$. Formally, we posit that there exists a constant $\ell>0$ such that 
\begin{itemize}
\item[(H3)] $\sup\nolimits_{x,t} \norm{\noisevol(x,t)}
	< \infty
	\quad
	\text{and}
	\quad
\norm{\noisevol(x',t) - \noisevol(x,t)}
	\leq \ell \norm{x'-x},$
\end{itemize}
where
\begin{equation}
\label{eq:Frobenius}
\norm{\noisevol}
	:= \sqrt{\tr\bracks{\noisevol\noisevol^{\top}}}
	= \sqrt{\sum\nolimits_{i=1}^{n} \sum\nolimits_{j=1}^{d} \abs{\noisevol_{ij}}^{2}}
\end{equation}
denotes the Frobenius (matrix) norm of $\noisevol$.
In particular, \eqref{eq:noise-reg} implies that there exists a positive constant $\noisevol_{\ast}\geq0$ such that
\begin{equation}
\label{eq:noisebound}
\norm{\noisevol(x,t)}^{2}
	\leq \noisevar
	\quad
	\text{for all $x\in\scrX$, $t\geq0$}.
\end{equation}
In what follows, it will be convenient to measure the intensity of the noise affecting \eqref{eq:SMD} via $\noisevol_{\ast}$;
of course, when $\noisevol_{\ast} = 0$, we recover the noiseless, deterministic dynamics \eqref{eq:MD}.

\section{Deterministic Analysis}
\label{sec:deterministic}
To establish a reference standard, we first focus on the deterministic regime of \eqref{eq:MD}, i.e. when $M(t)\equiv 0$ in \eqref{eq:SMD}. 

\subsection{Global Existence} 
We begin with a basic well-posedness result of \eqref{eq:MD}. 
%----------------------------------------------------------------------
\begin{proposition}
\label{prop:wp-det}
Under conditions \textnormal{(H1)} and \textnormal{(H2)}, the dynamical system \eqref{eq:MD} admits a unique solution from every initial condition $(s,y)\in\R_{+}\times\scrY$. 
\end{proposition}
%----------------------------------------------------------------------

%----------------------------------------------------------------------
{\it Proof.} 
Let $A(t,y) := -\lambda(t) v(Q(\eta(t)y))$ for all $t\in\Rp$, $y\in\scrY$.
Clearly, $A(t,y)$ is jointly continuous in $t$ and $y$.
Moreover, by (H2), $\lambda(t)$ has bounded first derivative and $\eta(t)$ is nonincreasing, so both $\lambda(t)$ and $\eta(t)$ are Lipschitz continuous.
Finally, by (H1), $v$ is $L$-Lipschitz continuous, implying in turn that
\begin{flalign}
%\label{eq:ALipschitz}
\dnorm{A(t,y_{1})-A(t,y_{2})}
	\leq\frac{L\eta(t)\lambda(t)}{\alpha} \norm{y_{1}-y_{2}}_{\ast}
	\quad
	\text{for all $y_{1},y_{2}\in\scrY$},
\end{flalign}
where $\alpha$ is the strong convexity constant of $h$, and we used \cref{prop:Q} to estimate the Lipschitz constant of $Q$. This shows that $A(t,y)$ is Lipschitz in $y$ for all $t$, so existence and uniqueness of local solutions follows from the Picard- Lindel\"{o}f theorem. (H2) further guarantees that the Lipschitz constant of $A(t,\cdot)$ can be chosen uniformly in $t$, so these solutions can be extended for all $t\geq0$.
\qed

%----------------------------------------------------------------------
Let $\bT:=\{(t,s)\vert 0\leq s\leq t\leq\infty\}$. Based on the above, we may define a \emph{non-autonomous semiflow} $Y:\bT\times\scrY \to\scrY$ satisfying (i) $Y(s,s,y)=y$ for all $s\geq 0$, (ii) $\frac{\partial Y(t,s,y)}{\partial t}=A(t,Y(t,s,y))$ for all $(t,s,y)\in\bT\times\scrY$, and (iii) $Y(t,s,Y(s,r,y))=Y(t,r,y)$ for $t\geq s\geq r\geq 0$. Since the dynamics will usually be started from an initial condition $(0,y)\in\R_{+}\times\scrY$, we will simplify the notation by writing $\phi(t,y)=Y(t,0,y)$ for all $(t,y)\in\R_{+}\times\scrY$. The resulting trajectory in the primal space is denoted by $\xi(t,y)=Q(\eta(t)\phi(t,y))$. Note that if $\lambda(t)$ and $\eta(t)$ are constant functions, the mapping $\phi(t,y)$ is the (autonomous) \emph{semiflow} of the dynamics \eqref{eq:MD}. 
%----------------------------------------------------------------------
\subsection{Convergence properties and performance}
Now, to analyze the convergence of \eqref{eq:MD}, we will consider two ``gap functions'' quantifying the distance between the primal trajectory, and the solution set of \eqref{eq:VI}:
%----------------------------------------------------------------------
\begin{itemize}
\item
In the general case, we will focus on the \emph{dual gap} function \cite{BorDut16}: 
\begin{equation}
\label{eq:gap}
g(x):= \max_{x'\in\scrX}\inner{v(x'),x - x'}.
\end{equation}
By \eqref{eq:Lip} and the compactness of $\scrX$, it follows that $g(x)$ is continuous, non-negative and convex;
moreover, we have $g(x) = 0$ if and only if $x$ is a solution of $\VI(\scrX,v)$ \cite[Proposition 3.1]{HarPan90}.
\item
For the saddle point problem \cref{ex:saddle}, we instead look at the \emph{Nikaido-Isoda gap function} \cite{NikIso55}:
\begin{equation}
\label{eq:gap-saddle}
G(p^{1},p^{2}):= \max_{x^{2}\in\scrX^{2}}U(p^{1},x^{2})
	- \min_{x^{1}\in\scrX^{1}} U(x^{1},p^{2}). 
\end{equation}
\end{itemize}
%----------------------------------------------------------------------
Since $U$ is convex-concave, it is immediate that $G(p^{1},p^{2})\geq g(p^{1},p^{2})$, where the operator involved in the definition of the dual gap function is given by the saddle-point operator \eqref{eq:v-saddle}. However, it is still true that $G(p^{1},p^{2})=0$ if and only if the pair $(p^{1},p^{2})$ is a saddle-point. Since both gap functions vanish only at solutions of \eqref{eq:VI}, we will prove trajectory convergence by monitoring the decrease of the relevant gap over time. This is achieved by introducing the so-called \emph{Fenchel coupling} \cite{Mer16}, an auxiliary energy function, defined as
\begin{equation}
\label{eq:Fench}
F(x,y):= h(x) + h^{\ast}(y) - \inner{y,x}
	\quad
	\text{for all $x\in\scrX$, $y\in\scrY$}.
\end{equation}
%where $h^{\ast}$ denotes the convex conjugate of $h$.
%\begin{remark}
%In a certain sense, the Fenchel coupling can be seen as a primal-dual extension of the well-known \emph{Bregman divergence} \cite{Bre67,Kiw97b}:
%\begin{equation}
%\label{eq:Breg}
%D(x,z):= h(x) - h(z) - h'(z;x-z)
%	\quad
%	\text{for all $x,z\in\scrX$}.
%\end{equation}
%More precisely, we have $F(x,y) \geq D(x,Q(y))$ with equality if and only if $Q(y)$ is interior \cite[Prop.~4.3]{Mer16}.
%\end{remark}
Some key properties of $F$ are summarized in the following proposition:
%----------------------------------------------------------------------
\begin{proposition}[\cite{Mer16}]
\label{prop:Fench}
Let $h$ be a distance-generating function on $\scrX$.
Then: 
\begin{enumerate}
[\itshape a\upshape)]
\item
$F(x,y)\geq \frac{\alpha}{2} \norm{Q(y)-x}^{2}$ for all $x\in\scrX$, $y\in\scrY$.
\item
Viewed as a function of $y$, $F(x,y)$ is convex, differentiable, and its gradient is given by
\begin{equation}
\nabla_{y} F(x,y)
	= Q(y) - x.
\end{equation}
\item
For al $x\in\scrX$ and all $y,y'\in\scrY$, we have
\begin{equation}
F(x,y')
	\leq F(x,y) + \inner{y'-y,Q(y)-x} + \frac{1}{2\alpha} \dnorm{y' - y}^{2}.
\end{equation}
\end{enumerate}
\end{proposition}
%----------------------------------------------------------------------

%We are now in a position to state and prove our first convergence result for \eqref{eq:MD}. 
In the sequel, if there is no danger of confusion, we will use the more concise notation $x(t)=\xi(t,y)$ and $y(t)=\phi(t,y)$, for the unique solution to \eqref{eq:MD} with initial condition $(0,y)\in\R_{+}\times\scrY$. Consider the averaged trajectory
\begin{flalign}
\label{eq:ergodic}
\bar x(t):= \frac{\int_{0}^{t} \lambda(s) x(s) \dif s}{\int_{0}^{t} \lambda(s) \dif s}= \frac{1}{S(t)} \int_{0}^{t} \lambda(s) x(s) \dif s,
\end{flalign}
where $S(t):= \int_{0}^{t} \lambda(s) \dif s.$ We then have the following convergence guarantee:

%----------------------------------------------------------------------
\begin{proposition}
\label{prop:gap-det}
Suppose that \eqref{eq:MD} is initialized at $(s,y)=(0,0)$, with resulting trajectories $y(t)=\phi(t,0)$ and $x(t)=\xi(t,0)$. Then:
\begin{equation}
\label{eq:gap-bound-det}
g(\bar x(t))
	\leq \frac{\scrD(h;\scrX)}{\eta(t)S(t)},
\end{equation}
where $\bar{x}(t)$ is the averaged trajectory constructed in \eqref{eq:ergodic}, and 
\begin{equation}
\label{eq:depth}
\scrD(h;\scrX):= \max_{x,x'\in\scrX} \{ h(x') - h(x) \}= \max h - \min h.
\end{equation}
In particular, if \eqref{eq:VI} is associated with a convex-concave saddle-point problem as in \cref{ex:saddle}, we have the guarantee:
\begin{equation}
\label{eq:saddle-bound-det}
G(\bar x(t))
	\leq \frac{\scrD(h_{1};\scrX^{1})+\scrD(h_{2};\scrX^{2})}{\eta(t)S(t)}.
\end{equation}
In both cases, whenever $\lim_{t\to\infty} \eta(t) S(t) = \infty$, $\bar x(t)$ converges to the solution set of $\VI(\scrX,v)$.
\end{proposition}
%----------------------------------------------------------------------

%----------------------------------------------------------------------
{\it Proof.} 
Given some $p\in\scrX$, let $H_{p}(t):= \frac{1}{\eta(t)}F(p,\eta(t)y(t))$. Then, with \cref{prop:Fench}, the fundamental theorem of calculus yields
\begin{equation}
\label{eq:H1}
H_{p}(t) - H_{p}(0)
	= -\int_{0}^{t}\lambda(s) \inner{v(x(s)),x(s)-p} \dif s
	- \int_{0}^{t} \frac{\dot{\eta}(s)}{\eta(s)^{2}}[h(p)-h(x(s))] \dif s,
\end{equation}
and, after rearranging, we obtain
\begin{flalign}
\int_{0}^{t}\lambda(s)\inner{v(x(s)),x(s)-p} \dif s
	&=H_{p}(0)-H_{p}(t)-\int_{0}^{t}\frac{\dot{\eta}(s)}{\eta(s)^{2}}[h(p)-h(x(s))]\dif s
	\notag\\
	&\leq  H_{p}(0)+\scrD(h;\scrX)\left(\frac{1}{\eta(t)}-\frac{1}{\eta(0)}\right).
\label{eq:reg-det}
\end{flalign}
Now, let $x_{c}:=\argmin\setdef{h(x)}{x\in\scrX}$ denote the ``prox-center'' of $\scrX$.
Since $\eta(0)>0$ and $y(0)=0$ by assumption, we readily get 
\begin{flalign}
H_{p}(0)
	= \frac{F(p,0)}{\eta(0)}
	=\frac{h(p) + h^{\ast}(0) - \inner{0,p}}{\eta(0)}
%	\notag\\
	=\frac{h(p)-h(x_{c})}{\eta(0)}
	\leq \frac{\scrD(h;\scrX)}{\eta(0)}.
\label{eq:H0}
\end{flalign}
From the monotonicity of $v$, we further deduce that
\begin{equation}
\label{eq:gap-det}
g(\bar x(t))
	\leq \frac{1}{S(t)} \max_{p\in\scrX} \int_{0}^{t} \lambda(s)\inner{v(x(s)),x(s)-p} \dif s.
%	= \frac{1}{S(t)}\delta(t). 
\end{equation}
Thus, substituting \eqref{eq:H0} in \eqref{eq:reg-det}, maximizing over $p\in\scrX$ and plugging the result into \eqref{eq:gap-det} gives \eqref{eq:gap-bound-det}.\\
Suppose now that \eqref{eq:VI} is associated to a convex-concave saddle-point problem, as in \cref{ex:saddle}.
In this case, we can replicate the above analysis for each component $x^{i}(t)$, $i=1,2$, of $x(t)$ to obtain the basic bounds
\begin{flalign}
\int_{0}^{t}\lambda(s)\inner{\nabla_{x^{1}}U(x(s)),x^{1}(s)-p^{1}} \dif s
	&\leq \frac{\scrD(h_{1};\scrX^{1})}{\eta(t)},
	\\
\int_{0}^{t}\lambda(s)\inner{-\nabla_{x^{2}}U(x(s)),x^{2}(s)-p^{2}} \dif s
	&\leq \frac{\scrD(h_{2};\scrX^{2})}{\eta(t)}.
\end{flalign}
Using the fact that $U$ is convex-concave, this leads to the value-based bounds
\begin{flalign}
\int_{0}^{t}\lambda(s)[U(x(s))-U(p^{1},x^{2}(s))] \dif s
	&\leq \frac{\scrD(h_{1};\scrX^{1})}{\eta(t)},
	\\
\int_{0}^{t}\lambda(s)[U(x^{1}(s),p^{2})-U(x(s))] \dif s
	&\leq \frac{\scrD(h_{2};\scrX^{2})}{\eta(t)}.
\end{flalign}
Summing these inequalities, dividing by $S(t)$, and using Jensen's inequality gives 
\begin{flalign}
U(\bar x^{1}(t),p^{2})-U(p^{1},\bar x^{2}(t))\leq\frac{\scrD(h_{1};\scrX^{1})+\scrD(h_{2};\scrX^{2})}{\eta(t)S(t)}
\end{flalign}
The bound \eqref{eq:saddle-bound-det} then follows by taking the supremum over $p^{1}$ and $p^{2}$, and using the definition of the Nikaido\textendash Isoda gap function.
\qed
%----------------------------------------------------------------------

The gap-based analysis of \cref{prop:gap-det} can be refined further in the case of \emph{strongly} monotone VI problems.

%----------------------------------------------------------------------
\begin{proposition}
Let $\sol$ denote the \textpar{necessarily unique} solution of a $\gamma$-strongly monotone problem $\VI(\scrX,v)$.
Then, with the same assumptions as in \cref{prop:gap-det}, we have
\begin{equation}
\label{eq:conv-strong-det}
\norm{\bar x(t) - \sol}^{2}
	\leq \frac{\scrD(h;\scrX)}{\gamma} \frac{1}{\eta(t)S(t)}.
\end{equation}
In particular, $\bar x(t)$ converges to $\sol$ whenever $\lim_{t\to\infty} \eta(t) S(t) = \infty$.
\end{proposition}
%----------------------------------------------------------------------

%----------------------------------------------------------------------
{\it Proof.} 
By Jensen's inequality, the strong monotonicity of $v$ and the assumption that $\sol$ solves $\VI(\scrX,v)$, we have:
\begin{flalign}
\gamma \norm{\bar x(t) - \sol}^{2}
	&\leq \frac{\gamma}{S(t)} \int_{0}^{t} \lambda(s) \norm*{x(s) - \sol}^{2} \dif s
	\tag{Jensen}
	\\
	&\leq \frac{1}{S(t)} \int_{0}^{t} \lambda(s) \inner{v(x(s)) - v(\sol), x(s) - \sol} \dif s
	\tag{$\gamma$-monotonicity}
	\\
	&\leq \frac{1}{S(t)} \int_{0}^{t} \lambda(s) \inner{v(x(s)), x(s) - \sol} \dif s
	\tag{optimality of $\sol$}
	\\
	&\leq \frac{\scrD(h;\scrX)}{\eta(t)S(t)},
\end{flalign}
where the last inequality follows as in the proof of \cref{prop:gap-det}.
The bound \eqref{eq:conv-strong-det} is then obtained by dividing both sides by $\gamma$.
\qed
%----------------------------------------------------------------------

The two results above are in the spirit of classical ergodic convergence results for monotone VI problems as in \cite{Bruck77,Nem04Prox,Nes09}. In particular, taking $\eta(t)=\sqrt{L/(2\alpha)}$ and $\lambda(t) = 1/(2\sqrt{t})$ gives the upper bound 
\begin{equation}
g(\bar x(t))
	\leq \scrD(h;\scrX) \sqrt{L/(\alpha t)},
\end{equation}
which is of the same order as the $\scrO(1/\sqrt{t})$ guarantees obtained in the references above.
However, the bound \eqref{eq:gap-det} does not have a term which is antagonistic to $\eta(t)$ or $\lambda(t)$, so, if \eqref{eq:MD} is run with constant $\lambda$ and $\eta$, we get an $\scrO(1/t)$ bound for $g(\bar x(t))$ (and/or $\norm{\bar x(t) - \sol}$ in the case of strongly monotone VI problems).%
\footnote{In fact, even faster convergence can be guaranteed if \eqref{eq:MD} is run with \emph{increasing} $\lambda(t)$.
In that case however, well-posedness is not immediately guaranteed, so we do not consider increasing $\lambda$ here.}
This suggests an important gap between continuous and discrete time;
for a similar phenomenon in the context of online convex optimization, see the regret minimization analysis of \cite{KM17}.

We close this section with a (nonergodic) trajectory convergence result for strictly monotone problems. For any path $X(\cdot):\R_{+}\to\scrX$, call the \emph{limit set} 
\begin{flalign}
\scrL\{X(\cdot)\}:=\bigcap_{t\geq 0}\cl[X([t,\infty)].
\end{flalign}
%----------------------------------------------------------------------
\begin{proposition}
\label{prop:conv-det}
Let $\sol$ denote the \textpar{necessarily unique} solution of a $\gamma$-strongly monotone problem $\VI(\scrX,v)$. Suppose Hypotheses \textnormal{(H1)} and \textnormal{(H2)} hold, and the parameters $\lambda$ and $\eta$ of \eqref{eq:MD} satisfy
\begin{flalign}
\textstyle
\inf_{t} \lambda(t) > 0
	\quad
	\text{and}
	\quad
\inf_{t} \eta(t) > 0.
\end{flalign}
Then, $\lim_{t\to\infty} \xi(t,y) = \sol$, for any $y\in\scrY$.
\end{proposition}
%----------------------------------------------------------------------

%----------------------------------------------------------------------
{\it Proof.} 
Let $x(t):=\xi(t,y)$ for $t\geq 0$, and assume that $\olim\in\scrL\{x(\cdot)\}$, but $\olim\neq\sol$.
%Since $\scrX_{\ast}$ is closed there exists an open neighborhood $O$ around $\olim$ such that $O\cap\scrX_{\ast}=\emptyset$.
Then, by assumption, there exists an open neighborhood $O$ of $\olim$ and a positive constant $a>0$ such that
\begin{equation}
\inner{v(x),x-\sol}
	\geq a
	\quad
	\text{for all $x\in O$}.
\end{equation}
Furthermore, since $\olim$ is an accumulation point of $x(t)$, there exists an increasing sequence $(t_{k})_{k\in\N}$ such that $t_{k}\uparrow\infty$ and $x(t_{k}) \to \olim$ as $k\to\infty$.
Thus, relabeling indices if necessary, we may assume without loss of generality that $x(t_{k})\in O$ for all $k\in\N$.
%Hence, there exists $M\in\N$ such that $x(t_{k})\in O$ for all $j\geq M$.

Now, for all $\eps>0$, we have 
\begin{flalign}
\norm{x(t_{k}+\eps)-x(t_{k})}
	&=\norm{Q(Y(t_{k}+\eps))-Q(Y(t_{k}))}
	\notag\\
	&\leq \frac{1}{\alpha}\norm{Y(t_{k}+\eps)-Y(t_{k})}_{\ast}
	\notag\\
	&\leq \frac{1}{\alpha}\int_{t_{k}}^{t_{k}+\eps}\lambda(s)\norm{v(x(s))}_{\ast}\dif s
	\notag\\
	&\leq \frac{1}{\alpha}\max_{x\in\scrX}\norm{v(x)}_{\ast}\int_{t_{k}}^{t_{k}+\eps}\lambda(s)\dif s
	\notag\\
	&\leq \frac{\eps\bar{\lambda}}{\alpha}\max_{x\in\scrX}\norm{v(x)}_{\ast},
\end{flalign}
where $\bar\lambda:= \lambda(0)$ denotes the maximum value of $\lambda(t)$. As this bound does not depend on $k$, we can choose $\eps>0$ small enough so that $x(t_{k}+s)\in O$ for all $s\in[0,\eps]$ and all $k\in\N$.
Thus, letting $H(t) := \eta(t)^{-1} F(\sol,\eta(t) y(t))$, and using \eqref{eq:reg-det}, we obtain
\begin{flalign}
%F(\sol,Y(t_{M+n})) - F(\sol,Y(t_{M}))
H(t_{n}) - H(t_{0})
	&= -\sum_{k=1}^{n} \int_{t_{k-1}}^{t_{k}} \lambda(s) \inner{v(x(s)),x(s) - \sol}\dif s
	+\scrD(h;\scrX)\left(\frac{1}{\eta(t_{n})}-\frac{1}{\eta(t_{0})}\right)
	\notag\\
	&\leq -a\underline{\lambda} \sum_{k=1}^{n}(t_{k}-t_{k-1})+\scrD(h;\scrX)\left(\frac{1}{\eta(t_{n})}-\frac{1}{\eta(t_{0})}\right)
	\notag\\
	&=-a\eps \underline{\lambda}n+\scrD(h;\scrX)\left(\frac{1}{\eta(t_{n})}-\frac{1}{\eta(t_{0})}\right),
\end{flalign}
where we have set $\underline\lambda:= \inf_{t}\lambda(t) > 0$. Given that $\inf_{t}\eta(t) > 0$, the above implies that $\lim_{n\to\infty} H(t_{n}) = -\infty$, contradicting the fact that $F(\sol,y)\geq 0$ for all $y\in\scrY$.
This implies that $\olim=\sol$; by compactness, $\scrL\{x(\cdot)\}\neq\emptyset$, so our claim follows.
\qed

\section{Analysis of the Stochastic Dynamics}
\label{sec:stochastic}
\subsection{Global Existence}
In this section, we turn to the stochastic system \eqref{eq:SMD}. As in the noise-free analysis of the previous section, we begin with a well-posedness result, stated for simplicity for deterministic initial conditions.

%----------------------------------------------------------------------
\begin{proposition}
\label{prop:wp-stoch}
Fix an initial condition $(s,y)\in\R_{+}\times \scrY$. Then, under Hypotheses \textnormal{(H1)-(H3)}, and up to a $\prob$-null set, the stochastic dynamics \eqref{eq:SMD} admit a unique strong solution $(Y(t))_{t\geq s}$ such that $Y(s) = y$.
\end{proposition}
%----------------------------------------------------------------------

%----------------------------------------------------------------------
{\it Proof.} 
Let $B(t,y):=- \lambda(t) \noisevol(Q(\eta(t)y),t)$ so \eqref{eq:SMD} can be written as
\begin{equation}
\label{eq:SDE}
dY(t)= A(t,Y(t))\dif t+B(t,Y(t)) \dif W(t),
\end{equation}
with $A(t,y)$ defined as in the proof of \cref{prop:wp-det}.
By (H2) and (H3), $B(t,y)$ inherits the boundedness and regularity properties of $\noisevol$;
In particular, (H2), (H3), together with Proposition \ref{prop:Q}c), imply that $B(t,y)$ is uniformly Lipschitz in $y$.
Under hypothesis (H1) and (H2), $A(t,y)$ is also uniformly Lipschitz in $y$ (cf. the proof of \cref{prop:wp-det}).
Our claim then follows by standard results in the well-posedness of stochastic differential equations \cite[Theorem~3.4]{Kha12}.
\qed
%----------------------------------------------------------------------
We denote by $Y(t,s,y)$ the unique strong solution of the It\^{o} stochastic differential equation \eqref{eq:SDE}, with initial condition $(s,y)\in\R_{+}\times \scrY$. As in the deterministic case, we are mostly interested in the process starting from the initial condition $(0,y)$, in which case we abuse the notation by writing $Y(t,y)=Y(t,0,y)$. The corresponding primal trajectories are generated by applying the mirror map $Q$ to the dual trajectories, so $X(t,y)=Q(\eta(t)Y(t,y))$, for all $(t,y)\in\R_{+}\times\scrY$. If there is no danger of confusion, we will consistently suppress the dependence on the initial position $y\in\scrY$ in both random processes. Clearly, if $\lambda(t)$ and $\eta(t)$ are constant function, the solutions of \eqref{eq:SMD} are time-autonomous.\\
We now give a brief overview of the results we obtain in this section. First, in \cref{sec:smallnoise}, we use the theory of \emph{asymptotic pseudo-trajectories} (APTs), developed by Benaïm and Hirsch \cite{BH96}, to establish almost sure trajectory convergence of \eqref{eq:SMD} to the solution of $\VI(\scrX,v)$, provided that $v$ is strictly monotone and the oracle noise in \eqref{eq:SMD} is vanishing at a rather slow, logarithmic rate. This strong convergence result relies heavily on the shadowing property of the dual trajectory and its deterministic counterpart $\phi(t,y)$. (see Section \ref{sec:smallnoise}).
On the other hand, if the driving noise process is persistent, we cannot expect the primal trajectory $X(t)$ to converge - some averaging has to be done in this case.
Thus, following a long tradition on ergodic convergence for mirror descent, we investigate in \cref{sec:ergodic} the asymptotics of a weighted time-average of $X(t)$.
Finally, we complement our ergodic convergence results with a large deviation principle showing that the ergodic average of $X(t)$ is exponentially concentrated around its mean (\cref{sec:LDP}).

%----------------------------------------------------------------------
%% SMALL NOISE
%----------------------------------------------------------------------
\subsection{The Small Noise Limit}
\label{sec:smallnoise}
We begin with the case where the oracle noise in \eqref{eq:SMD} satisfies the asymptotic decay condition $\norm{\noisevol(x,t)} \leq \beta(t)$ for some nonincreasing function $\beta\from\Rp\to\Rp$ such that
\begin{itemize}
\item[(H4)]
 $\int_{0}^{\infty} \exp\parens[\bigg]{-\frac{c}{\beta^{2}(t)}} \dif t < \infty
	\quad
	\text{for all $c>0$}$. 
\end{itemize}
For instance, this condition is trivially satisfied if $\noisevol(x,t)$ vanishes at a logarithmic rate, i.e. $\beta(t) = o(1/\sqrt{\log(t)})$.
For technical reasons, we will also need the additional ``Fenchel reciprocity'' condition
\begin{itemize}
\item[(H5)]
$F(p,y_{n})\to 0
	\quad
	\text{whenever}
	\quad
	Q(y_{n})\to p$.
\end{itemize}
Under the decay rate condition \textnormal{(H4)},
%and the reciprocity requirement \textnormal{(H5)},
and working for simplicity with constant $\eta(t) = \lambda(t) = 1$, results of \cite[Proposition 4.1]{BH96} imply that any strong solution $Y(t)$ of \eqref{eq:SMD} is an (APT) of the deterministic dynamics \eqref{eq:MD} in the following sense:
%----------------------------------------------------------------------
\begin{definition}
\label{def:APT}
Assume that $\eta(t) = \lambda(t) = 1$, for all $t\geq 0$. Let $\phi\from\Rp\times\scrY\to\scrY$, $(t,y)\mapsto\phi(t,y)$, denote the semiflow induced by \eqref{eq:MD} on $\scrY$. A continuous curve $Y\from\Rp\to\scrY$ is said to be an \emph{asymptotic pseudo-trajectory} (APT) of \eqref{eq:MD}, if
\begin{equation}
\label{eq:APT}
\tag{APT}
\lim_{t\to\infty} \sup_{0\leq s\leq T} \dnorm{Y(t+s) - \phi(s,Y(t))}
	= 0
	\quad
	\text{for all $T>0$.}
\end{equation}
\end{definition}
%----------------------------------------------------------------------

In words, \cref{def:APT} states that an APT of \eqref{eq:MD} tracks the solutions of \eqref{eq:MD} to arbitrary accuracy over arbitrarily long time windows. Thanks to this property, we are able to establish the following global convergence theorem for \eqref{eq:SMD} with vanishing oracle noise:

\begin{theorem}
\label{thm:smallnoise}
Assume that $v$ is strictly monotone, and let $\sol$ denote the \textpar{necessarily unique} solution of $\VI(\scrX,v)$. If Hypotheses \textnormal{(H1)-(H5)} hold, and \eqref{eq:SMD} is run with $\lambda(t) = \eta(t) = 1$, we have
\begin{equation}
\probof*{\lim_{t\to\infty} \norm{X(t,y) - \sol} = 0}
	=1
	\quad\forall y\in\scrY.
\end{equation}
\end{theorem}
%----------------------------------------------------

The proof of \cref{thm:smallnoise} requires some auxiliary results, which we provide below. We begin with a strong recurrence result for neighborhoods of the (unique) solution $\sol$ of $\VI(\scrX,v)$ under \eqref{eq:MD}:

%----------------------------------------------------------------------
\begin{lemma}
\label{lem:recurrentdeterministic}
With assumptions as in \cref{thm:smallnoise}, let $\nhd$ be an open neighborhood of $\sol$ in $\scrX$ and let $\xi(t,y)=Q(\eta(t)\phi(t,y))$. Define the stopping time
\begin{equation}
t_{\nhd}(y)
	:= \inf\setdef{t\geq0}{\xi(t,y)\in \nhd}.
\end{equation}
Then $t_{\nhd}(y)<\infty$ for all $y\in\scrY$. 
\end{lemma}
%----------------------------------------------------------------------

%----------------------------------------------------------------------
{\it Proof.} 
%Fix the initial condition $y\in\scrY$ arbitrary.
%To simplify notation we just write $x(t)=\xi(t,y)$ and $y(t)=\phi(t,y)$. 
%Define the mapping $H(t):=F(\sol,y(t))$. By the chain-rule of calculus 
Fix the initialization $y\in\scrY$ of \eqref{eq:MD}, and let $y(t):= \phi(t,y),\;x(t):=Q(\phi(t,y))$ denote the induced solutions of \eqref{eq:MD}, and let $H(t) := F(\sol,y(t))$.
Then, by \cref{prop:Fench}, and the chain rule applied to \eqref{eq:MD}, we get 
\begin{equation}
H(t)
	= H(0) - \int_{0}^{t} \inner{v(x(s)),x(s)-\sol} \dif s. 
\end{equation}
Since $v$ is strictly monotone and $\sol$ solves $\VI(\scrX,v)$, there exists some $a\equiv a_{\nhd} > 0$ such that
\begin{equation}
\inner{v(x),x-\sol}
	\geq a
	\quad
	\text{for all $x\in\scrX\setminus \nhd$}.
\end{equation}
Hence, if $t_{\nhd}(y) = \infty$, we would have
\begin{equation}
H(t)
	\leq H(0) - a t
	\quad
	\text{for all $t\geq0$},
\end{equation}
implying in turn that $\lim_{t\to\infty} H(t) = -\infty$. This contradicts the fact that $H(t)\geq0$, so we conclude that $t_{\nhd}(y) < \infty$.
\qed
%----------------------------------------------------------------------
Next, we extend this result to the stochastic regime:
%----------------------------------------------------------------------
\begin{lemma}
\label{lem:recurrentnoise}
With assumptions as in \cref{thm:smallnoise}, let $\nhd$ be an open neighborhood of $\sol$ in $\scrX$ and define the stopping time
\begin{equation}
\tau_{\nhd}(y):=\inf\{t\geq 0: X(t,y)\in \nhd\},
\end{equation}
Then, $\tau_{\nhd}(y)$ is almost surely finite for all $y\in\scrY$.
\end{lemma}
%----------------------------------------------------------------------

%----------------------------------------------------------------------
{\it Proof.} 
Suppose there exists some initial condition $y_{0}\in\scrY$ such that $\probof{\tau_{\nhd}(y_{0})=\infty}>0$.
Then there exists a measurable set $\Omega_{0}\subseteq\Omega$, with $\probof{\Omega_{0}}>0$, and such that $\tau_{\nhd}(\omega,y_{0})=\infty$ for all $\omega\in\Omega_{0}$. Now, define $H(t):=F(\sol,Y(t,y_{0}))$ and set $X(t)=X(t,y_{0})$. By the weak It\^{o} lemma \eqref{eq:Ito} proven in \cref{app:stoch}, we get 
\begin{equation}
H(t) - H(0)
	\leq -\int_{0}^{t} \inner{v(X(s)),X(s) - \sol}\dif s
	+ \frac{1}{2\alpha} \int_{0}^{t} \norm{\noisevol(X(s),s)}^{2}\dif s
	+ I_{\sol}(t)
\end{equation}
where $I_{\sol}(t):= \int_{0}^{t} \inner{X(s) - \sol,\noisevol(X(s))\cdot\dif W(s)}$ is a continuous local martingale.
Since $v$ is strictly monotone, the same reasoning as in the proof of \cref{lem:recurrentdeterministic} yields
\begin{equation}
H(t)
	\leq H(0)
	- at
	+ I_{\sol}(t)
	+ \frac{\noisevar}{2\alpha}
\end{equation}
for some $a \equiv a_{\nhd} > 0$ and for all $t\in[0,\tau_{\nhd}(y))$. Furthermore, by an argument based on the law of the iterated logarithm and the Dambis-Dubins-Schwarz time-change theorem for martingales as in the proof of \cref{thm:gap-stoch}, we get
\begin{equation}
\text{$I_{\sol}(t)/t \to 0$ almost surely as $t\to\infty$}.
\end{equation}
Combining this with the estimate for $H(t)$ above, we get $\lim_{t\to\infty} H(t) = -\infty$ for $\prob$-almost all $\omega\in\Omega_{0}$. This contradicts the fact that $H(t)\geq0$, and our claim follows.
\qed

The above result shows that the primal process $X(t)$ hits any neighborhood of $\sol$ in finite time \as.
Thanks to this important recurrence property, we are finally in a position to prove \cref{thm:smallnoise}:

{\it Proof of Theorem \ref{thm:smallnoise}.}
Fix some $\eps>0$ and let $N_{\eps}:=\setdef{x=Q(y)}{F(\sol,y)<\eps}$. Let $y\in\scrY$ be arbitrary.
We first claim, that there exists a deterministic time $T\equiv T(\eps)$ such that
$F(\sol,\phi(T,y))\leq\max\{\eps,F(\sol,y)+\eps\}$. Indeed, consider the hitting time 
\begin{equation}
t_{\eps}(y):=\inf\setdef{t\geq 0}{x(t)\in N_{\eps}},
\end{equation}
where $x(t) :=Q(\phi(t,y))$. By \textnormal{(H5)}, $N_{\eps}$ contains a neighborhood of $\sol$;
hence, by \cref{lem:recurrentdeterministic}, we have $t_{\eps}(y) < \infty$.
Moreover, observe that
\begin{equation}
\label{eq:Fdecrease}
\frac{d}{dt} F(\sol,\phi(t,y))
	=-\inner{v(x(t)),x(t) - \sol}
	\leq 0
	\quad
	\text{for all $y\in\scrY$}.
\end{equation}
% For a given $y\in\scrY$ and $t\in[0,t_{\eps}(y))$,
The strict monotonicity of $v$ and the fact that $\sol$ solves \eqref{eq:VI} implies that there exists a positive constant $\kappa \equiv \kappa_{\eps} >0$ such that $\inner{v(x),x-\sol} \geq \kappa$ for all $x\in\scrX\setminus N_{\eps}$.
 Hence, combining this with \eqref{eq:Fdecrease}, we readily see that
\begin{equation}
F(\sol,\phi(t,y)) - F(\sol,y)
	\leq -\kappa t
	\quad
	\text{for all $t\in[0,t_{\eps}(y))$}.
\end{equation}
Now, set $T = \eps/\kappa$.
If $T<t_{\eps}(y)$, we immediately conclude that  
\begin{equation}
F(\sol,\phi(T,y)) - F(\sol,y)
	\leq-\eps.
\end{equation}
Otherwise, if $T \geq t_{\eps}(y)$, we again use the descent property \eqref{eq:Fdecrease} to get 
\begin{equation}
F(\sol,\phi(T,y))
	\leq F(\sol,\phi(t_{\eps}(y),y))
	\leq\eps. 
\end{equation}
In both cases we have $F(\sol,\phi(T,y)) \leq \max\{\eps,F(\sol,y)-\eps\}$, as claimed.

To proceed, pick $\delta\equiv \delta_{\eps}>0$ such that 
\begin{equation}
\label{eq:delta}
\delta_{\eps}\diam(\scrX)+\frac{\delta^{2}_{\eps}}{2\alpha}<\eps,
\end{equation}
where $\diam(\scrX):= \max\setdef{\norm{x'-x}_{2}}{x,x'\in\scrX}$ denotes the Euclidean diameter of $\scrX$.
By Proposition 4.1 of \cite{BH96}, the strong solution $Y$ of \eqref{eq:SDE} (viewed as a stochastic flow) is an APT of the deterministic semiflow $\phi$ with probability $1$.
Hence, we can choose an \as finite random time $\theta_{\eps}$ such that $\sup_{s\in[0,T]}\dnorm{Y(t+s)-\phi(s,Y(t))}\leq\delta_{\eps}$ for all $t\geq\theta_{\eps}$.
Combining this with item (c) of \cref{prop:Fench}, we then get
%for all $t\geq \theta_{\eps}(y),\;s\in[0,T]$ 
\begin{flalign}
F(\sol,Y(t+s,y))
	&\leq F(\sol,\phi(s,Y(t,y)))
	\notag\\
	&+\inner{Y(t+s,y)-\phi(s,Y(t,y)),Q(\phi(s,Y(t,y)))-\sol}
	\notag\\
	&+\frac{1}{2\alpha}\dnorm{Y(t+s,y)-\phi(s,Y(t,y))}^{2}
	\notag\\
	&\leq F(\sol,\phi(s,Y(t,y)))+\delta_{\eps}\diam(\scrX)+\frac{\delta_{\eps}^{2}}{2\alpha}
	\notag\\
	&\leq F(\sol,\phi(s,Y(t,y)))+\eps,
\end{flalign}
where the last inequality follows from the estimate \eqref{eq:delta}.

Now, choose a random time $T_{0}\geq\max\{\theta_{\eps}(y),t_{\eps}(y)\}$ and $T=\eps/\kappa$ as above.
Then, by definition, we have $F(\sol,Y(T_{0},y))\leq2\eps$ with probability $1$.
Hence, for all $s\in[0,T]$, we get
\begin{equation}
F(\sol,Y(T_{0}+s,y))
	\leq F(\sol,\phi(s,Y(T_{0},y))) + \eps\leq F(\sol,Y(T_{0},y)) + \eps\leq 3\eps. 
\end{equation}
Since $F(\sol,\phi(T,Y(T_{0},y))) \leq \max\{\eps,F(\sol,Y(T_{0},y)) - \eps\}\leq \eps$, we also get 
\begin{flalign}
F(\sol,Y(T_{0}+T+s,y))
	&\leq F(\sol,\phi(s,Y(T_{0}+T,y))) + \eps
	\notag\\
	&\leq F(\sol,Y(T_{0}+T,y)) + \eps
	\notag\\
	&\leq 3\eps,
\end{flalign}
and hence
\begin{equation}
F(\sol,Y(T_{0}+s,y))\leq 3\eps
	\quad
	\text{for all $s\in[T,2T]$}.
\end{equation}
Using this as the basis for an induction argument, we readily get 
\begin{equation}
F(\sol,Y(T_{0}+s,y))
	\leq 3\eps
	\quad
	\text{for all $s\in[nT,(n+1)T]$},
\end{equation}
with probability $1$.
Since $\eps$ was arbitrary, we obtain $F(\sol,Y(t,y))\to 0$, implying in turn that $X(t)\to\sol$ \as by \cref{prop:Fench}.
\qed
%----------------------------------------------------------------------

%----------------------------------------------------------------------
%% ERGODIC
%----------------------------------------------------------------------
\subsection{Ergodic Convergence}
\label{sec:ergodic}

We now proceed with an ergodic convergence result in the spirit of \cref{prop:gap-det}. The results presented in this section are derived under the assumption that \eqref{eq:SMD} is started with the initial conditions $(s,y)=(0,0)$. This is only done to make the presentation clearer; see Remark \ref{rem:init}.\\
Set $S(t):= \int_{0}^{t} \lambda(s) \dif s$, and $L(t):= \sqrt{\int_{0}^{t} \lambda^{2}(s) \dif s}$,
and let
\begin{equation}
\label{eq:ergodic-stoch}
\bar X(t):= \frac{1}{S(t)}\int_{0}^{t}\lambda(s)X(s)\dif s,
\end{equation}
denote the ``ergodic average'' of $X(t)=X(t,0,0)$.
%----------------------------------------------------------------------
\begin{theorem}
\label{thm:gap-stoch}
%With notation and assumptions as in \cref{lem:gap-stoch}, we have:
Under Hypotheses \textnormal{(H1)-(H3)}, we have:
\begin{equation}
\label{eq:gap-rate-stoch}
g(\bar X(t))
	= \bigoh\parens*{\frac{1}{\eta(t) S(t)}}
	+ \bigoh\parens*{\frac{\int_{0}^{t} \lambda^{2}(s) \eta(s)\dif s}{S(t)}}
	+ \bigoh\parens*{\frac{L(t) \sqrt{\log\log L(t)}}{S(t)}},
\end{equation}
with probability $1$.
In particular, $\bar X(t)$ converges \as to the solution set of $\VI(\scrX,v)$ provided that
\begin{inparaenum}
[\itshape a\upshape)]
\item
$\lim_{t\to\infty} \eta(t) S(t) = \infty$;
and
\item
$\lim_{t\to\infty} \eta(t) \lambda(t) = 0$.
\end{inparaenum}
\end{theorem}
%----------------------------------------------------------------------

Before discussing the proof of \cref{thm:gap-stoch}, it is worth noting the interplay between the two variable weight parameters, $\lambda(t)$ and $\eta(t)$.
In particular, if \eqref{eq:SMD} is run with weight sequences of the form $1/t^{q}$ for some $q>0$, we obtain:

%----------------------------------------------------------------------
\begin{corollary}
\label{cor:gap-stoch}
Suppose that \eqref{eq:SMD} is run with $\lambda(t) = (1+t)^{-a}$ and $\eta(t) = (1+t)^{-b}$ for some $a,b\in[0,1]$.
Then, with assumptions as in \cref{thm:gap-stoch}, we have:
\begin{equation}
\label{eq:gap-rate-explicit}
g(\bar X(t))
	= \tilde\bigoh\parens*{t^{a+b-1}}
	+ \tilde\bigoh\parens*{t^{-a-b}}
	+ \tilde\bigoh\parens*{t^{-1/2}}.
\end{equation}
\end{corollary}
%----------------------------------------------------------------------

In the above, the $\tilde\bigoh(\cdot)$ notation signifies ``$\bigoh(\cdot)$ up to logarithmic factors''.%
\footnote{More precisely, we write $f(x) = \tilde\bigoh(g(x))$ when $f(x) = \bigoh(g(x) \log^{k}g(x))$ for some $k>0$.}
Up to such factors, \eqref{eq:gap-rate-explicit} is optimized when $a+b=1/2$;
if these factors are to be considered, any choice with $a+b=1/2$ and $b>0$ gives the same rate of convergence, indicating that the role of the post-multiplication factor $\eta(t)$ is crucial to finetune the convergence rate of \eqref{eq:SMD}.
We find this observation particularly appealing as it is reminiscent of Nesterov's remark that ``running the discrete-time algorithm \eqref{eq:MD-discrete} with the \emph{best} step-size strategy $\lambda_{t}$ and fixed $\eta$ [\dots] gives the same (infinite) constant as the corresponding strategy for fixed $\lambda$ and variable $\eta_{t}$'' \cite[p.~224]{Nes09}.

The proof of \cref{thm:gap-stoch} relies crucially on the following lemma, which provides an explicit estimate for the decay rate of the employed gap functions.
%----------------------------------------------------------------------
\begin{lemma}
\label{lem:gap-stoch}
If \eqref{eq:SMD} is initialized at $(0,0)$, and Hypotheses \textnormal{(H1)-(H3)} hold, then:
\begin{equation}
\label{eq:gap-stoch}
g(\bar X(t))
	\leq \frac{\scrD(h;\scrX)}{\eta(t)S(t)}
	+ \frac{\noisevar}{2\alpha} \frac{\int_{0}^{t} \lambda^{2}(s) \eta(s)\dif s}{S(t)}
	+\frac{I(t)}{S(t)} 
\end{equation}
where
$I(t):= \sup_{p\in\scrX} I_{p}(t)$ and 
\begin{equation}
\label{eq:Ip}
I_{p}(t):= \int_{0}^{t} \lambda(s) \inner{p - X(s), \noisevol(X(s),s) \cdot dW(s)}.
\end{equation}
If \eqref{eq:VI} is associated with a convex-concave saddle-point problem as in \cref{ex:saddle}, we have
\begin{equation}
\label{eq:saddle-bound-stoch}
G(\bar{X}^{1}(t),\bar{X}^{2}(t))
	\leq \frac{\scrD_{\textup{sp}}}{\eta(t)S(t)}
	+\frac{\noisevar}{2\alpha_{\textup{sp}}} \frac{\int_{0}^{t} \lambda^{2}(s) \eta(s)\dif s}{S(t)}
	+\frac{J(t)}{S(t)},
\end{equation}
where we have set
$\scrD_{\textup{sp}}:= \scrD(h_{1};\scrX^{1})+\scrD(h_{2};\scrX^{2})$,
$1/\alpha_{\textup{sp}}:= 1/\alpha_{1} +1/\alpha_{2}$,
and
$J(t):= \sup_{p^{1}\in\scrX^{1},p^{2}\in\scrX^{2}} \{I_{p^{1}}(t) + I_{p^{2}}(t)\}$.
\end{lemma}
%----------------------------------------------------------------------

\begin{remark}
\label{rem:init}
The initialization assumption in \cref{lem:gap-stoch} is not crucial:
we only make it to simplify the explicit expression \eqref{eq:gap-stoch}.
If \eqref{eq:SMD} is initialized at a different point, the proof of \cref{lem:gap-stoch} shows that the bound \eqref{eq:gap-stoch} is correct only up to $\bigoh(1/S(t))$.
Since all terms in \eqref{eq:gap-stoch} are no faster than $\bigoh(1/S(t))$, initialization plays no role in the proof of \cref{thm:gap-stoch} below.
\end{remark}

%----------------------------------------------------------------------
{\it Proof of \cref{lem:gap-stoch}.}
Fix some $p\in\scrX$ and let $H_{p}(t):= \eta(t)^{-1}F(p,\eta(t)Y(t))$ as in the proof of \cref{prop:gap-det}.
Then, by the weak It\^{o} formula \eqref{eq:Ito} in \cref{app:stoch}, we have
\begin{flalign}
H_{p}(t)
	\leq H_{p}(0)
	&-\int_{0}^{t} \frac{\dot{\eta}(s)}{\eta(s)^{2}} H_{p}(s)\dif s
	+\frac{1}{\eta(t)} \int_{0}^{t} \inner{X(s)-p,\dot\eta(s) Y(s)} \dif s
	\notag\\
	&+ \int_{0}^{t} \inner{X(s) - p,\dif Y(s)}
	+ \frac{1}{2\alpha} \int_{0}^{t} \lambda^{2}(s) \eta(s)\norm{\noisevol(X(s),s)}^{2} \dif s.
\label{eq:Hp-stoch-1}
\end{flalign}
To proceed, let
\begin{equation}
\label{eq:Rp}
R_{p}(t):= \int_{0}^{t}\lambda(s) \inner{v(X(s)),X(s)-p} \dif s,
\end{equation}
so 
\begin{flalign}
\int_{0}^{t} \inner{X(s)-p, \dif Y(s)} &= -\int_{0}^{t} \lambda(s) \inner{X(s)-p,v(X(s)) \dif s + \dif M(s)}
	\notag\\
	&=  -R_{p}(t) + I_{p}(t), 
\end{flalign}
with $I_{p}(t)$ given by \eqref{eq:Ip}.
Then, rearranging and bounding the second term of \eqref{eq:Hp-stoch-1} as in the proof of \cref{prop:gap-det}, we obtain
\begin{flalign}
R_{p}(t)
	&\leq H_{p}(0)-H_{p}(t)
	+ \scrD(h,\scrX) \left(\frac{1}{\eta(t)}-\frac{1}{\eta(0)}\right)
	\notag\\
	&+ I_{p}(t)
	+ \frac{1}{2\alpha} \int_{0}^{t} \lambda^{2}(s) \eta(s) \norm{\noisevol(X(s),s)}^{2} \dif s
	\notag\\
	&\leq H_{p}(0)
	+ \scrD(h;\scrX) \left(\frac{1}{\eta(t)} - \frac{1}{\eta(0)}\right)
	+ I_{p}(t)
	+ \frac{\noisevar}{2\alpha} \int_{0}^{t} \lambda^{2}(s) \eta(s) \dif s.
\end{flalign}
With \eqref{eq:SMD} initialized at $y=0$, \cref{eq:H0} gives $H_{p}(0) \leq \scrD(h;\scrX)/\eta(0)$.
Thus, by Jensen's inequality and the monotonicity of $v$, we get
\begin{flalign}
\inner{v(p),\bar X(t) - p}
	&= \frac{1}{S(t)} \int_{0}^{t} \lambda(s) \inner{v(p),X(s)-p} \dif s
	\notag\\
	&\leq \frac{1}{S(t)} \int_{0}^{t} \lambda(s) \inner{v(X(s)),X(s)-p} \dif s
	= \frac{R_{p}(t)}{S(t)}
	\notag\\
	&\leq \frac{\scrD(h;\scrX)}{\eta(t)S(t)}
	+ \frac{\noisevar}{2\alpha} \frac{\int_{0}^{t} \lambda^{2}(s) \eta(s)\dif s}{S(t)}
	+\frac{I_{p}(t)}{S(t)}.
\label{eq:bound1}
\end{flalign}
The bound \eqref{eq:gap-stoch} then follows by noting that $g(\bar X(t)) = \max_{p\in\scrX} \inner{v(p),\bar{X}(t)-p}$.

Now, assume that \eqref{eq:VI} is associated to a convex-concave saddle point problem as in \cref{ex:saddle}.
As in the proof of \cref{prop:gap-det}, we first replicate the analysis above for each component of the problem, and we then sum the two components to get an overall bound for the Nikaido-Isoda gap function $G$.
Specifically, applying \eqref{eq:bound1} to \eqref{eq:v-saddle}, we readily get
\begin{equation}
\label{eq:bound1-saddle}
\int_{0}^{t} \lambda(s) \inner{v^{i}(X(s)),X^{i}(s) - p^{i}}
	\leq \frac{\scrD(h_{i};\scrX^{i})}{\eta(t) S(t)}
	+ \frac{\noisevar}{2\alpha^{i}} \frac{\int_{0}^{t} \lambda^{2}(s) \eta(s)\dif s}{S(t)}
	+\frac{I_{p^{i}}(t)}{S(t)},
\end{equation}
where $i\in\{1,2\}$.
Moreover, Jensen's inequality yields
\begin{flalign}
U(\bar X^{1}(t),p^{2}) - U(p^{1},\bar X^{2}(t))
	&\leq \frac{1}{S(t)} \int_{0}^{t} \lambda(s) \bracks{U(X^{1}(s),p^{2}) - U(p^{1},X^{2}(s))} \dif s
	\notag\\
	&\leq \frac{1}{S(t)} \int_{0}^{t} \lambda(s) \inner{\nabla_{x^{1}}U(X(s)),X^{1}(s) - p^{1}} \dif s
	\notag\\
	&- \frac{1}{S(t)} \int_{0}^{t} \lambda(s) \inner{\nabla_{x^{2}}U(X(s)),X^{2}(s) - p^{2}} \dif s
	\notag\\
	&\leq \frac{\scrD_{\textup{sp}}}{\eta(t)S(t)}
	+\frac{\noisevar}{2\alpha_{\textup{sp}}} \frac{\int_{0}^{t} \lambda^{2}(s) \eta(s)\dif s}{S(t)}
	+\frac{I_{p^{1}}(t) + I_{p^{2}}(t)}{S(t)},
\end{flalign}
with the last inequality following from \eqref{eq:bound1-saddle}. Our claim then follows by maximizing over $(p^{1},p^{2})$ and recalling the definition \eqref{eq:gap-saddle} of the Nikaido-Isoda gap function.
\qed
%----------------------------------------------------------------------

%With \cref{lem:gap-stoch} at hand, we proceed to establish the convergence of \eqref{eq:SMD} in the ergodic sense.
Clearly, the crucial unknown in the bound \eqref{eq:gap-stoch} is the stochastic term $I(t)$. To obtain convergence of $\bar X(t)$ to the solution set of $\VI(\scrX,v)$, the term $I(t)$ must grow slower than $S(t)$.
As we show below, this is indeed the case:

%----------------------------------------------------------------------
{\it Proof of \cref{thm:gap-stoch}.}
By \cref{lem:gap-stoch} and \cref{rem:init}, it suffices to show that the term $I(t)$ grows as $\bigoh(L(t)\sqrt{\log\log L(t)})$ with probability $1$.
To do so, let $\kappa_{p} := \bracks{I_{p}}$ denote the quadratic variation of $I_{p}$.%
\footnote{Recall here that the quadratic variation of a stochastic process $M(t)$ is the continuous increasing and progressively measurable process, defined as $\bracks{M(t)} = \lim_{\abs{\Pi}\to0} \sum_{1\leq j \leq k} (M(t_{j}) - M(t_{j-1}))^{2}$, where the limit is taken over all partitions $\Pi = \{t_{0} = 0 < t_{1} < \dotsi < t_{k} = t\}$ of $[0,t]$ with mesh $\abs{\Pi} \equiv \max_{j} \abs{t_{j} - t_{j-1}} \to 0$ \cite{KS98}.}
Then, the rules of stochastic calculus yield
\begin{flalign}
d\kappa_{p}(t)
	&= dI_{p}(t) \cdot dI_{p}(t)
	\notag\\
	&= \lambda^{2}(t) \sum_{i,j=1}^{n} \sum_{k=1}^{d} (X_{i}(t) - p_{i}) (X_{j}(t) - p_{j}) \noisevol_{ik}(X(t),t) \noisevol_{jk}(X(t),t) \dif t
	\notag\\
	&\leq \norm{X(t) - p}_{2}^{2} \noisevar \lambda^{2}(t)
	\notag\\
	&\leq \diam(\scrX)^{2} \noisevar \lambda^{2}(t),
\end{flalign}
where $\diam(\scrX):= \max\setdef{\norm{x'-x}_{2}}{x,x'\in\scrX}$ denotes the Euclidean diameter of $\scrX$.
Hence, for all $t\geq0$, we get the quadratic variation bound
\begin{equation}
\label{eq:cov-bound}
\kappa_{p}(t)
	\leq \diam(\scrX)^{2} \noisevar \int_{0}^{t} \lambda^{2}(s) \dif s
	= \bigoh(L^{2}(t)).
\end{equation}

Now, let $\kappa_{p}(\infty) := \lim_{t\to\infty} \kappa_{p}(t) \in[0,\infty]$ and set
\begin{equation}
\tau_{p}(s): = \begin{cases}
	\inf\setdef{t\geq0}{\kappa_{p}(t) > s}
		&\quad
		\text{if $s\leq\kappa_{p}(\infty)$},
		\\
	\infty
		&\quad
		\text{otherwise}.
	\end{cases}
\end{equation}
The process $\tau_{p}(s)$ is finite, non-negative, non-decreasing and right-continuous on $[0,\kappa_{p}(\infty))$;
moreover, it is easy to check that $\kappa_{p}(\tau_{p}(s)) = s \wedge \kappa_{p}(\infty)$ and $\tau_{p}(\kappa_{p}(t)) = t$ \cite[Problem~3.4.5]{KS98}.
Therefore, by the Dambis-Dubins-Schwarz time-change theorem for martingales \cite[Theorem.~3.4.6 and Problem.~3.4.7]{KS98}, there exists a standard, one-dimensional Wiener process $(B_{p}(t))_{t\geq0}$ adapted to a modified filtration $\tilde\scrF_{s} = \scrF_{\tau_{p}(s)}$ (possibly defined on an extended probability space), and such that $B_{p}(\kappa_{p}(t)) = I_{p}(t)$ for all $t\geq0$ (except possibly on a $\prob$-null set). Hence, for all $t>0$, we have
\begin{equation}
\frac{I_{p}(t)}{S(t)}
	= \frac{B_{p}(\kappa_{p}(t))}{S(t)}
	= \frac{B_{p}(\kappa_{p}(t))}{\sqrt{\kappa_{p}(t) \log\log\kappa_{p}(t)}}
	\times \frac{\sqrt{\kappa_{p}(t) \log\log\kappa_{p}(t)}}{S(t)}.
\end{equation}
By the law of the iterated logarithm \cite{KS98}, the first factor above is bounded almost surely;
as for the second, \eqref{eq:cov-bound} gives $\sqrt{\kappa_{p}(t) \log\log\kappa_{p}(t)} = \bigoh(L(t) \sqrt{\log\log L(t)})$.
Thus, combining all of the above, we get
\begin{equation}
\frac{I(t)}{S(t)}
	= \frac{\max_{p\in\scrX} I_{p}(t)}{S(t)}
	= \bigoh\parens*{\frac{L(t) \sqrt{\log\log L(t)}}{S(t)}},
\end{equation}
so our claim follows from \eqref{eq:gap-stoch}.

To complete our proof, note first that the condition $\lim_{t\to\infty} \eta(t) S(t) = \infty$ implies that $\lim_{t\to\infty} S(t) = \infty$ (given that $\eta(t)$ is nonincreasing).
Thus, by de l'H\^{o}pital's rule and the assumption $\lim_{t\to\infty} \lambda(t) \eta(t) = 0$, we also get $S(t)^{-1} \int_{0}^{t} \lambda^{2}(s) \eta(s) \dif s = 0$.
Finally, for the last term of \eqref{eq:gap-rate-stoch}, consider the following two cases:
\begin{enumerate}
\item
If $\lim_{t\to\infty} L(t) < \infty$, we trivially have $\lim_{t\to\infty} L(t) \sqrt{\log\log L(t)} \big/ S(t) = 0$ as well.
\item
Otherwise, if $\lim_{t\to\infty} L(t) = \infty$, de l'H\^{o}pital's rule readily yields
\begin{equation}
\lim_{t\to\infty} \frac{L^{2}(t)}{S^{2}(t)}
	= \lim_{t\to\infty} \frac{\lambda^{2}(t)}{2 \lambda(t) S(t)}
	= \frac{1}{2} \lim_{t\to\infty} \frac{\lambda(t)}{S(t)}
	= 0,
\end{equation}
by the boundedness of $\lambda(t)$.
Another application of de l'H\^{o}pital's rule gives
\begin{equation}
\lim_{t\to\infty} \frac{L^{3}(t)}{S^{2}(t)}
	= \lim_{t\to\infty} \frac{(L^{2}(t))^{3/2}}{S^{2}(t)}
	= \frac{3}{4} \lim_{t\to\infty} \frac{\lambda^{2}(t) L(t)}{\lambda(t) S(t)}
	= \frac{3}{4} \lim_{t\to\infty} \frac{\lambda(t) L(t)}{S(t)}
	= 0,
\end{equation}
so
\begin{equation}
\limsup_{t\to\infty} \frac{L(t)\sqrt{\log\log L(t)}}{S(t)}
	\leq \limsup_{t\to\infty} \sqrt{\frac{L^{3}(t)}{S^{2}(t)}}
	= 0.
\end{equation}
\end{enumerate}
The above shows that, under the stated assumptions, the RHS of \eqref{eq:gap-rate-stoch} converges to $0$ almost surely, implying in turn that $\bar X(t)$ converges to the solution set of $\VI(\scrX,v)$ with probability $1$.
\qed
%----------------------------------------------------------------------

%----------------------------------------------------------------------
%% LARGE DEVIATIONS
%----------------------------------------------------------------------
\subsection{Large deviations}
\label{sec:LDP}

In this section we study the concentration properties of \eqref{eq:SMD} in terms of the dual gap function. As in the previous section, we will assume that \eqref{eq:SMD} is issued from the initial condition $(s,y)=(0,0)$. \\
First, recall that for every $p\in\scrX$ we have the upper bound 
\begin{equation}
R_{p}(t)
	\leq \frac{\scrD(h;\scrX)}{\eta(t)}
	+\frac{\noisevar}{2\alpha}\int_{0}^{t}\lambda^{2}(s)\eta(s)\dif s
	+I_{p}(t).
\end{equation}
with $R_{p}(t)$ and $I_{p}(t)$ defined as in \eqref{eq:Rp} and \eqref{eq:Ip} respectively.
Since $I_{p}(t)$ is a continuous martingale starting at $0$, we have $\exof{I_{p}(t)} = 0$, implying in turn that
%Hence, for every $p\in\scrX$, using the fact that $\inner{v(p),\bar{X}(t)-p}\leq \frac{R_{p}(t)}{S(t)}$, we readily get 
\begin{equation}
\exof{\inner{v(p),\bar{X}(t)-p}}
	\leq \frac{\scrD(h;\scrX)}{S(t)\eta(t)}
	+ \frac{\noisevar}{2\alpha S(t)}\int_{0}^{t}\lambda^{2}(s)\eta(s)\dif s
%	+ I_{p}(t)
	=\frac{K(t)}{2S(t)},
\end{equation}
where 
\begin{equation}
\label{eq:K}
K(t):= \frac{2\scrD(h;\scrX)}{\eta(t)}+ \frac{\noisevar}{\alpha} \int_{0}^{t}\lambda^{2}(s)\eta(s)\dif s.
\end{equation}
%Hence, taking the supremum over all $p\in\scrX$ and using Jensen's inequality, we get the mean gap bound 
%\begin{equation}
%\exof{g(\bar{X}(t)}
	%\leq \frac{K(t)}{2S(t)}.
%\end{equation}
Markov's inequality therefore implies that
\begin{equation}
\label{eq:Markov}
\probof{\inner{v(p),\bar{X}(t)-p} \geq \delta}
	\leq\frac{1}{\delta} \frac{K(t)}{2S(t)}
	\quad
	\text{for all $\delta>0$}.
\end{equation}

The bound \eqref{eq:Markov} provides a first estimate of the probability of observing a large gap from the solution of \eqref{eq:VI}, but because it relies only on Markov's inequality, it is rather crude.
To refine it, we provide below a ``large deviations'' bound that shows that the ergodic gap process $g(\bar X(t))$ is exponentially concentrated around its mean value:

%----------------------------------------------------------------------
\begin{theorem}
\label{thm:LDP}
Suppose \textnormal{(H1)-(H3)} hold, and that \eqref{eq:SMD} is started from the initial condition $(s,y)=(0,0)$. Then, for all $\delta>0$ and all $t>0$, we have 
\begin{equation}
\label{eq:gap-LDP}
\probof{g(\bar{X}(t))\geq \scrQ_{0}(t)+\delta\scrQ_{1}(t)}
	\leq \exp(-\delta^{2}/4),
\end{equation}
where 
\begin{subequations}
\begin{flalign}
\label{eq:Q0}
\scrQ_{0}(t):=\frac{K(t)}{S(t)},
\intertext{and}
\label{eq:Q1}
\scrQ_{1}(t):=\frac{\sqrt{\kappa}\noisedev\diam(\scrX)L(t)}{S(t)},
\end{flalign}
\end{subequations}
with $\kappa>0$ a positive constant depending only on $\scrX$ and $\norm{\cdot}$.
\end{theorem}

The concentration bound \eqref{eq:gap-LDP} can also be formulated as follows:
\begin{corollary}
With notation and assumptions as in \cref{thm:LDP}, we have
\begin{flalign}
g(\bar X(t))
	&\leq \scrQ_{0}(t) + 2 \scrQ_{1}(t) \sqrt{\log(1/\delta)}
	\notag\\
	&= \bigoh\parens*{\frac{1}{\eta(t) S(t)}}
	+ \bigoh\parens*{\frac{\int_{0}^{t} \lambda^{2}(s) \eta(s)\dif s}{S(t)}}
	+ \bigoh\parens*{\frac{L(t) \sqrt{\log(1/\delta)}}{S(t)}},
\label{eq:gap-LDP-highprob}
\end{flalign}
with probability at least $1-\delta$.
In particular, if \eqref{eq:SMD} is run with parameters $\lambda(t) = (1+t)^{-a}$ and $\eta(t) = (1+t)^{-b}$ for some $a,b\in[0,1]$, we have
\begin{equation}
\label{eq:gap-LDP-explicit}
g(\bar X(t))
	= \bigoh\parens*{t^{a+b-1}}
	+ \bigoh\parens*{t^{-a-b}}
	+ \bigoh\parens*{t^{-1/2}},
\end{equation}
with arbitrarily high probability.
\end{corollary}

%----------------------------------------------------------------------
To prove \cref{thm:LDP}, define first the auxiliary processes 
\begin{flalign}
Z(t):= \int_{0}^{t}\lambda(s)\noisevol(X(s),s)\dif W(s)
	\quad
	\text{and}
	\quad
	P(t):= Q(\eta(t)Z(t)).
\end{flalign}
We then have:

%----------------------------------------------------------------------
\begin{lemma}
\label{lem:aux1}
For all $p\in\scrX$ we have 
\begin{equation}
\label{eq:UB1}
\int_{0}^{t} \lambda(s)\inner{p-P(s),\noisevol(X(s),s)\dif W(s)}
	\leq \frac{\scrD(h;\scrX)}{\eta(t)}
	+ \frac{\noisevar}{2\alpha}\int_{0}^{t}\lambda^{2}(s)\eta(s)\dif s.
\end{equation}
\end{lemma}
%----------------------------------------------------------------------

%----------------------------------------------------------------------
{\it Proof.} 
The proof follows the same lines as Lemma \ref{lem:gap-stoch}.
Specifcially, given a reference point $p\in\scrX$, define the process $\tilde H_{p}(t):= \frac{1}{\eta(t)} F(p,\eta(t)Z(t))$.
Then, by the weak It\^{o} formula \eqref{eq:Ito} in \cref{app:stoch}, we have 
\begin{flalign}
\tilde H_{p}(t)
	&\leq \tilde H_{p}(0)
	\notag\\
	&-\int_{0}^{t} \frac{\dot{\eta}(s)}{\eta(s)^{2}}\tilde H_{p}(s)\dif s+\frac{1}{\eta(t)} \int_{0}^{t} \inner{\xi(s)-p,\dot\eta(s)Z(s)} \dif s
	\notag\\
	&+ \int_{0}^{t} \inner{P(s) - p,\dif Z(s)}
	+ \frac{1}{2\alpha} \int_{0}^{t} \lambda^{2}(s) \eta(s)\norm{\noisevol(X(s),s)}^{2} \dif s
	\notag\\
	&\leq -\int_{0}^{t}\frac{\dot{\eta}(s)}{\eta(s)}[h(p)-h(P(s))]\dif s
	\notag\\
	&+ \int_{0}^{t}\lambda(s)\inner{P(s)-p,\noisevol(X(s),s)\dif W(s)}
	+ \frac{\noisevar}{2\alpha} \int_{0}^{t} \lambda^{2}(s)\eta(s) \dif s. 
\end{flalign}
We thus get
\begin{flalign}
\int_{0}^{t}\lambda(s)\inner{p-P(s),\noisevol(X(s),s)\dif W(s)}
	&\leq \tilde H_{p}(0)
	+ \scrD(h;\scrX) \, \parens*{\frac{1}{\eta(t)}-\frac{1}{\eta(0)}}
	\notag\\
	&+ \int_{0}^{t}\frac{\lambda^{2}(s)\eta(s)}{2\alpha}\noisevar \dif s
	\notag\\
	&\leq \frac{\scrD(h;\scrX)}{\eta(t)}
	+ \int_{0}^{t}\frac{\lambda^{2}(s)\eta(s)}{2\alpha}\noisevar\dif s,
\end{flalign} 
as claimed.
\qed
%----------------------------------------------------------------------

We are now ready to establish our large deviations principle for \eqref{eq:SMD}:

{\it Proof of Theorem \ref{thm:LDP}.}
For $p\in\scrX$ and $t>0$ fixed, we have 
\begin{flalign}
R_{p}(t)
	&\leq \frac{\scrD(h;\scrX)}{\eta(t)}
	+ \frac{\noisevar}{2\alpha}\int_{0}^{t}\lambda^{2}(s)\eta(s)\dif s
	+ \int_{0}^{t}\lambda(s)\inner{p-X(s),\noisevol(X(s),s)\dif W(s)}
	\notag\\
	&= \frac{\scrD(h;\scrX)}{\eta(t)}
	+ \frac{\noisevar}{2\alpha}\int_{0}^{t}\lambda^{2}(s)\eta(s)\dif s
	+ \int_{0}^{t}\lambda(s)\inner{p-P(s),\noisevol(X(s),s)\dif W(s)}
	\notag\\
	&+ \int_{0}^{t}\lambda(s)\inner{P(s)-X(s),\noisevol(X(s),s)\dif W(s)}
	\notag\\
	&\leq \frac{2\scrD(h;\scrX)}{\eta(t)}
	+ \frac{\noisevar}{\alpha}\int_{0}^{t}\lambda^{2}(s)\eta(s)\dif s
	+ \int_{0}^{t}\inner{P(s)-X(s),\noisevol(X(s),s)\dif W(s)},
\end{flalign}
where we used \eqref{eq:UB1} to obtain the last inequality.
To proceed, let 
\begin{equation}
\Delta(t):=\int_{0}^{t} \lambda(s)\inner{P(s)-X(s),\noisevol(X(s),s) \dif W(s)}.
\end{equation}
The process $\Delta(t)$ is a continuous martingale starting at $0$ which is almost surely bounded in $L^{2}$, thus providing an upper bound for $R_{p}(t)$ which is independent of the reference point $p\in\scrX$.
Indeed, recalling the definition \eqref{eq:K} of $K(t)$,
% $K(t)=2\scrD(h;\scrX)/\eta(t) + \alpha^{-1}\noisevar \int_{0}^{t} \lambda^{2}(s)\eta(s)\dif s$,
 we see that 
\begin{equation}
R_{p}(t)\leq K(t)+\Delta(t),
\end{equation}
so
\begin{flalign}
g(\bar{X}(t))
	\leq \frac{K(t) +\Delta(t)}{S(t)}
	\quad
	\text{for all $t>0$}.
\end{flalign}
In turn, this implies that
for all $\eps,t>0$,
%\begin{equation}
%\{g(\bar{X}(t)\geq\eps\}\subseteq \{\Delta(t)\geq \eps S(t)-K(t)\},
%\end{equation}
%and hence 
\begin{equation}
\probof{g(\bar{X}(t))\geq\eps}
	\leq \probof{\Delta(t)\geq\eps S(t)-K(t)}
	\quad
	\text{for all $\eps,t>0$}.
\end{equation}

To prove the theorem, we are left to bound the right-hand side of the above expression.
To that end, letting $\rho(t):= [\Delta(t),\Delta(t)]$ denote the quadratic variation of $\Delta(t)$, the Cauchy-Schwarz inequality readily gives 
\begin{flalign}
\exof{\exp(\theta \Delta(t))}
	&= \exof{\exp(\theta\Delta(t)-b\rho(t))\exp(b\rho(t))}
	\notag\\
	&\leq \sqrt{\exof{\exp(2\theta\Delta(t)-2b\rho(t))}}
		\sqrt{\exof{\exp(2b\rho(t))}}.
\end{flalign}
Setting $b=\theta^{2}$, the expressions inside the first expected value is just the stochastic exponential of the process $2\theta\Delta(t)$.
Moreover, a straightforward calculation shows that 
\begin{flalign}
\rho(t)
	&= \int_{0}^{t} \lambda^{2}(s) \norm{\noisevol(X(s),s) \cdot (P(s)-X(s))}_{2}^{2} \dif s
	\notag\\
	&\leq \kappa \int_{0}^{t}\lambda^{2}(s)\norm{\noisevol(X(s),s)}^{2}\norm{P(s)-X(s)}^{2}\dif s
	\notag\\
	&\leq\kappa\noisevar\diam(\scrX)^{2}L^{2}(t),
\end{flalign}
where $\kappa>0$ is a universal constant accounting for the equivalence of the Euclidean norm $\norm{\cdot}_{2}$ and the primal norm $\norm{\cdot}$ on $\scrX$.
The above implies that $\rho(t)$ is bounded over every compact interval, showing that Novikov's condition is satisfied (see e.g. \cite{KS98}). We conclude that the process $\exp(2\theta\Delta(t)-2\theta^{2}\rho(t))$ is a true martingale with expected value $1$.
Hence, letting $\varphi(t):=\kappa\noisevar\diam(\scrX)^{2}L^{2}(t)$, we get
\begin{equation}
\label{eq:bounddelta}
\exof{\exp(\theta\Delta(t))}
	\leq \sqrt{\exof{\exp(2\theta^{2}\rho(t)}}
	\leq \exp(\theta^{2}\varphi(t)).
\end{equation}
Combining all of the above, we see that for all $a>0$
\begin{flalign}
\probof{\Delta(t)\geq a}
	&= \probof{\exp(\theta\Delta(t))\geq \exp(\theta a)}
	\notag\\
	&\leq \exp(-a\theta) \exof{\exp(\theta\Delta(t))}
	\tag{Markov inequality}
	\\
	&=\exp(-a\theta + \theta^{2} \varphi(t))
	\quad
	\text{for all $a>0$},
\end{flalign}
with the last line following from \eqref{eq:bounddelta}.
Minimizing the above with respect to $\theta$ then gives
\begin{equation}
\probof{\Delta(t)\geq a}
	\leq \exp\left(-\frac{a^{2}}{4\varphi(t)}\right).
\end{equation}
Hence, by unrolling \cref{eq:Q0,eq:Q1}, we finally obtain the bound
\begin{flalign}
\probof{g(\bar{X}(t))\geq \scrQ_{0}(t)+\delta\scrQ_{1}(t)}
	&\leq\probof{\Delta(t)\geq \scrQ_{0}(t)S(t) + \delta\scrQ_{1}(t)S(t)-K(t)}
	\notag\\
	&\leq\probof{\Delta(t)\geq\delta \sqrt{\varphi(t)}}
	\notag\\
	&\leq \exp(-\delta^{2}/4),
\end{flalign}
for all $\delta>0$, as claimed.
\qed

%\begin{remark}
%The concentration bound established in Theorem \ref{thm:LDP} can also be formulated in a slightly different, though equivalent, way. A quick inspection shows that the Theorem is saying that for every $\delta>0$ and $t>0$, with probability at least $1-\delta$, the dual gap function is bounded by 
%\begin{flalign}
%g(\bar{X}(t))\leq \scrQ_{0}(t)+2\scrQ_{1}(t)\sqrt{\ln(1/\delta)}.
%\end{flalign}
%\end{remark}

\section{Conclusion}
This paper examined a continuous-time dynamical system for solving monotone variational inequality problems with random inputs. The key element of our analysis is the identification of a energy-type function, which allows us to prove ergodic convergence of generated trajectories in the deterministic as well as in the stochastic case. Future research should extend the present work in the following dimensions. First, it is not clear yet how the continuous-time method will help us in the derivation of a consistent numerical scheme. A naive Euler-discretization might potentially lead to a loss in speed of convergence (see \cite{WibWilJor16}). Second, it is of great interest to relax the monotonicity assumption we made on the involved operator. We are currently investigating these extensions. Third, it is of interest to consider different noise models as well. In particular, it would be interesting to know how the results derived in this paper change when the stochastic perturbation come from a jump Markov process, or more generally, a L\'{e}vy process. This extension would likely need new techniques and we regard this as an important contribution for future work. 
 
\section*{acknowledgements}
We thank the co-Editor in Chief, Professor Franco Giannessi, for clarifications on the conceptual differences between Minty and Stampacchia Variational inequality problems, and some pointers to the relevant literature. 
P.~Mertikopoulos was partially supported by
the French National Research Agency (ANR) grant
ORACLESS (ANR\textendash 16\textendash CE33\textendash 0004\textendash 01)
and
the COST Action CA16228 "European Network for Game Theory" (GAMENET).
The research of M. Staudigl is partially supported by the COST Action CA16228 "European Network for Game Theory" (GAMENET).

\appendix  %This command ends the counting of sections.
\section*{Appendix A: Results from Convex Analysis}\label{app:computation}
In this appendix we collect some simple facts on the analysis of convex differentiable functions with Lipschitz continuous gradients. Denote by $\textbf{C}^{1,1}_{L}(\Rn)$ the totality of such functions, with $L$ being the Lipschitz constant of the gradient mapping $\nabla\psi$:
\begin{flalign}
\norm{\psi(y+\delta)-\psi(y)}_{\ast}\leq L\norm{\delta}_{\ast}\qquad\forall y,\delta\in\Rn. 
\end{flalign}

\begin{proposition}\label{prop:Hessianbound}
Let $\psi\in\textbf{C}^{1,1}_{L}(\Rn)$ be convex. Then $\psi$ is almost everywhere twice differentiable with Hessian $\nabla^{2}\psi$  and  
\begin{equation}
0\leq \nabla^{2}\psi(y)\leq L\Id.\quad\Leb-\text{a.e.}.
\end{equation}
\end{proposition}

{\it Proof.}
For every $\psi\in\textbf{C}^{1,1}_{L}(\Rn)$, the well-known descent lemma (\cite{Nes04}, Theorem 2.1.5) implies that
\begin{equation}\label{eq:descent}
\psi(y+\delta)\leq \psi(y)+\inner{\nabla\psi(y),\delta}+\frac{L}{2}\norm{\delta}_{\ast}^{2} \qquad\forall y,\delta\in\Rn. 
\end{equation}
By Alexandrov's theorem (see, e.g., \cite{YonZho99}, Lemma 6.6), it follows that $\psi$ is $\Leb$-almost everywhere twice differentiable. Hence, there exists a measurable set $\Lambda$ such that $\Leb(\Lambda)=0$, and for all $\bar{y}\in\Rn\setminus\Lambda$ there exists $(p,P)\in\Rn\times\R^{n\times n}_{sym}$ such that 
\begin{equation}\label{eq:psi1}
\psi(\bar{y}+y)=\psi(\bar{y})+\inner{p,y}+\frac{1}{2}\inner{Py,y}+\theta(\bar{y},y),
\end{equation}
where $\lim_{\norm{y}_{\ast}\to 0}\frac{\theta(\bar{y},y)}{\norm{y}^{2}_{\ast}}=0$. We have $p=\nabla\psi(\bar{y})$ and identify $P$ with the a.e. defined Hessian $\nabla^{2}\psi(\bar{y})$. On the other hand, convexity implies
\begin{equation}\label{eq:psi2}
\psi(\bar{y}+y)\geq \psi(\bar{y})+\inner{\nabla \psi(\bar{y}),y}\qquad\forall \bar{y},y\in\Rn. 
\end{equation}
Choosing $y=t e$, where $e\in\Rn$ is an arbitrary $\norm{\cdot}_{\ast}$-unit vector and $t>0$, it follows 
\begin{flalign}
-\frac{1}{t^{2}}\theta(\bar{y},te)&\stackrel{\eqref{eq:psi2}}\leq \frac{1}{t^{2}}[\psi(\bar{y}+te)-\psi(\bar{y})-\inner{\nabla\psi(\bar{y}),te}]-\frac{1}{t^{2}}\theta(\bar{y},te)
	\notag\\
	&\stackrel{\eqref{eq:psi1}}=\frac{1}{2}\inner{\nabla^{2}\psi(\bar{y})e,e}
	\notag\\
	&\stackrel{\eqref{eq:descent}}\leq  \frac{L}{2}-\frac{1}{t^{2}}\theta(\bar{y},te).
\end{flalign}
Letting $t\to 0^{+}$ we get
\begin{equation}\label{eq:estimateh}
0\leq \frac{1}{2}\inner{\nabla^{2}\psi(\bar{y})e,e}\leq \frac{L}{2}\qquad\forall \bar{y}\in\Rn\setminus\Lambda,
\end{equation}
which implies $\nabla^{2}\psi(\bar{y})\leq L\Id$. 
\qed

\section*{Appendix B: Results from Stochastic Analysis}
\label{app:stoch}
The following result is the generalized It\^{o} formula used in the main text. 
\begin{proposition}
\label{prop:ito}
Let $Y$ be an It\^{o} process in $\Rn$ of the form 
\begin{equation}
Y_{t}=Y_{0}+\int_{0}^{t}F_{s}\dif s+\int_{0}^{t}G_{s}\dif W(s).
\end{equation}
Let $\psi\in\textbf{C}^{1,1}_{L}(\Rn)$ be convex. Then, for all $t\geq 0$, we have 
\begin{flalign}
\psi(Y_{t})\leq\psi(Y_{0})+\int_{0}^{t}\inner{\nabla\psi(Y_{s}),\dif Y_{s}}+\frac{L}{2}\int_{0}^{t}\norm{G_{s}}^{2}\dif s
\end{flalign}
\end{proposition}
{\it Proof.}
Since $\psi\in\textbf{C}^{1,1}_{L}(\Rn)$ is convex, Proposition \ref{prop:Hessianbound} shows that $\psi$ is almost everywhere twice differentiable with Hessian $\nabla^{2}\psi$. Furthermore, this Hessian matrix satisfies $0\leq\nabla^{2}\psi(y)\leq L\Id$, for all $y\in\Rn$ outside a set of Lebesgue measure 0.\\
Introduce the mollifier 
\begin{flalign}
\rho(u):=\left\{\begin{array}{lr} c\exp\left(\frac{-1}{1-\norm{u}_{\ast}^{2}}\right) & \text{ if }\norm{u}_{\ast}<1,\\
0 & \text{if }\norm{u}_{\ast}\geq 1.
\end{array}\right. 
\end{flalign}
Choose the constant $c>0$ so that $\int_{\Rn}\rho(u)\dif u=1$. For every $\eps>0$ define 
\begin{flalign}
&\rho_{\eps}(u):=\eps^{-n}\rho(u/\eps),\\
&\psi_{\eps}(y):=\psi\circledast\rho_{\eps}(y):=\int_{\Rn}\psi(y-u)\rho_{\eps}(u)\dif u.
\end{flalign}
Then, $\psi_{\eps}\in\textbf{C}^{\infty}(\Rn)$ and the standard form of It\^{o}'s formula gives us
\begin{flalign}
\psi_{\eps}(Y_{t})&=\psi_{\eps}(Y_{s})+\int_{s}^{t}\inner{\nabla\psi_{\eps}(Y_{r}),\dif Y_{r}}+\frac{1}{2}\int_{s}^{t}\tr\left[\nabla^{2}\psi_{\eps}(Y_{r})G_{r}G_{r}^{\top}\right]\dif r
	\notag\\
	&=\psi_{\eps}(Y_{s})+\int_{s}^{t}\inner{\int_{\Rn}\nabla\psi(z)\rho_{\eps}(Y_{r}-z)\dif z,\dif Y_{r}}
	\notag\\
	&+\frac{1}{2}\int_{s}^{t}\int_{\Rn}\tr\left[\nabla^{2}\psi(z)G_{r}G_{r}^{\top}\right]\rho_{\eps}(Y_{r}-z)\dif r\dif z.
\end{flalign}
Since $\tr(\nabla^{2}\psi(z)G_{r}G_{r}^{\top})\leq L\norm{G_{r}}^{2}$, we get 
\begin{flalign}
\psi_{\eps}(Y_{t})&\leq \psi_{\eps}(Y_{s})+\int_{s}^{t}\inner{\int_{\Rn}\nabla\psi(z)\rho_{\eps}(Y_{r}-z)\dif z,\dif Y_{r}}+\frac{L}{2}\int_{s}^{t}\norm{G_{r}}^{2}\dif r. 
\end{flalign}
Letting $\eps\to 0^{+}$, using the uniform convergence of the involved data, proves the result.
\qed

Applying this result to the dual process of \eqref{eq:SMD} and using \eqref{eq:estimateh}, gives for $F_{s}=A(s,Y(s,y))$ and $G_{s}=B(s,Y(s,y))$, the following version of the generalized It\^{o} rule: 
\begin{equation}
\label{eq:Ito}
\psi(Y_{t})\leq\psi(Y_{0})+\int_{0}^{t}\inner{\nabla\psi(Y_{s}),\dif Y_{s}}\dif s+\frac{1}{2\alpha}\int_{0}^{t}\norm{B(s,Y(s,y))}^{2}\dif s
\end{equation}

%*************************************************************
%*****    BIBLIOGRAPHY
%*************************************************************
\bibliographystyle{plain}

\end{document}